\documentclass{amsart}
\usepackage{amscd,amssymb,amsxtra}
\usepackage[mathscr]{eucal}
\usepackage[all]{xy}

\usepackage[usenames,dvipsnames]{color}

\definecolor{softyellow}{rgb}{.8745,.8588,.7647}
\definecolor{softblue}{rgb}{.6862,.8078,.8745}
\definecolor{softgreen}{rgb}{.8000,.8745,.8156}

\usepackage{graphicx}
\usepackage{mathtools}

\usepackage{hyperref}
\hypersetup{%
  bookmarksnumbered=true,%
  bookmarks=true,%
  colorlinks=true,%
 allcolors=black,%
}

\title{Real Wilson Spaces I}

\author{M.~A.~Hill}
\address{Department of Mathematics \\ UCLA
\\ Los Angeles, CA}
\email{mikehill@math.ucla.edu}

\author{M.~J.~Hopkins}
\address{Department of Mathematics \\ Harvard University
\\Cambridge, MA 02138}
\email{mjh@math.harvard.edu}

\newtheorem{thm}[equation]{Theorem}
\newtheorem{cor}[equation]{Corollary}
\newtheorem{lem}[equation]{Lemma}
\newtheorem{prop}[equation]{Proposition}
\newtheorem*{thm*}{Theorem}
\newtheorem*{cor*}{Corollary}
\newtheorem*{lem*}{Lemma}
\newtheorem*{prop*}{Proposition}

\ifx\undefined\theoremstyle
	\newtheorem{definition}[equation]{Definition}
	\newenvironment{defin}{\begin{definition}\rm}{\end{definition}}
	\newtheorem{conj}[equation]{Conjecture}
	\newtheorem{rem}[equation]{Remark}
	\newtheorem*{rem*}{Remark}
	\newtheorem{rems}[equation]{Remarks}

\else
	\theoremstyle{definition}

	\newtheorem{defin}[equation]{Definition}
	
	\newtheorem{notation}[equation]{Notation}

	\theoremstyle{remark}

	\newtheorem{rem}[equation]{Remark}
	\newtheorem*{rem*}{Remark}

\ifx\undefined\prem	\newtheorem{prem}{\bf Working Remark} \fi

\fi
\newtheorem{eg}[equation]{Example}

\ifx\undefined\pf
\newenvironment{pf}{\bigskip{\em Proof:\/}}{\qed\medskip}
\newenvironment{pf*}[1]{\bigskip{\em #1:\/}}{\qed\medskip}
\fi

\ifx\undefined\numberwithin
\def\numberwithin#1#2{\makeatletter\@ifundefined{c@#1}{\@nocnterrr}{%
  \@ifundefined{c@#2}{\@nocnterr}{%
  \@addtoreset{#1}{#2}%
  \toks@\expandafter\expandafter\expandafter{\csname the#1\endcsname}%
  \expandafter\xdef\csname the#1\endcsname
    {\expandafter\noexpand\csname the#2\endcsname
     .\the\toks@}}}\makeatother}\fi

\DeclareMathOperator{\tor}{tor}

\newcommand{\Q}{{\mathbb Q}}
\newcommand{\Z}{{\mathbb Z}}
\newcommand{\R}{{\mathbb R}}

\newcommand{\zerowidth}[1]{\hbox to 0pt{\hss$\displaystyle #1$\hss}}
\newcommand{\LL}{[\mkern-2mu[}
\newcommand{\RR}{]\mkern-2mu]}

\newcommand{\C}{\mathbb C}

\ifx\undefined\eqref\newcommand{\eqref}[1]{\rm (\ref{#1})}\fi

\newcounter{thmItem}

\newcommand{\thmItemref}[1]
 	 {{\rm \ref{#1})}}

\newenvironment{thmList}{\begin{list}%
{\rm \roman{thmItem})}{\usecounter{thmItem}
\setlength{\labelwidth}{2em}
\setlength{\itemindent}{2em}
\setlength{\leftmargin}{0pt}
\setlength{\listparindent}{0pt}
\setlength{\parsep}{0pt}
\setlength{\partopsep}{0pt}
\setlength{\itemsep}{\medskipamount}
\setlength{\topsep}{\medskipamount}
}}{\end{list}}

\newcounter{textItem}

\newcounter{condItem}

\ifx\undefined\slot
\newcommand{\slot}{\,-\,}
\fi

\newcommand{\MU}{MU}
\newcommand{\BU}{BU}
\newcommand{\umu}{\underline{\MU}}
\newcommand{\mur}{{\MU_{\R}}}
\newcommand{\bur}{{\BU_{\R}}}
\newcommand{\umur}{\underline{\mur}}

\newcommand{\hringh}{\mathscr{H}}
\newcommand{\hringr}{\mathscr{R}}

\newcommand{\hringmu}{\mathscr{MU}}
\newcommand{\hringmur}{\mathscr{MUR}}
\newcommand{\hringmurw}{\hringmu_{\text{RW}}}
\newcommand{\rwHringIdeal}{K_{\text{RW}}}

\newcommand{\freeHopfRing}{\mathscr{F}}

\newcommand{\rbullet}[1]{{#1}^{\bullet}}

\newcommand{\zt}{C_{2}}
\newcommand{\ft}{\mathbb{F}_{2}}
\newcommand{\ftwo}{\ft}
\newcommand{\hft}{H\ft}

\newcommand{\Zm}{\underline{\Z}}
\newcommand{\zm}{\Zm}

\DeclareMathOperator{\coalg}{\mathbf{coalg}}

\DeclareMathOperator{\hopfalg}{\mathbf{HopfAlg}}

\DeclareMathOperator{\hopfrings}{\mathbf{HopfRings}}

\newcommand{\hopfRingsMU}{\hopfrings_{\rbullet{MU}}}

\newcommand{\cp}{\mathbf{CP}}
\newcommand{\rp}{\mathbf{RP}}

\DeclareMathOperator{\sym}{Sym}
\newcommand{\A}{\mathbb{A}}

\newcommand{\gfp}[1]{\Phi^{#1}}
\newcommand{\bgfp}[1]{{\bar{\Phi}^{#1}}}
\newcommand{\phig}{\gfp{\zt}}
\newcommand{\bphig}{\bgfp{\zt}}

\DeclareMathOperator{\spaces}{\mathfrak{Spaces}}

\DeclareMathOperator{\ho}{ho}

\newcommand{\rmod}[1]{\modules_{#1}}

\newcommand{\zmod}{\modules_{H\Z}}
\newcommand{\zmodmap}[2]{\zmod(#1,#2)}
\newcommand{\hzmod}{\modules_{H\Zm}}
\newcommand{\gmodules}[1]{\modules^{\ast}_{#1}}

\newcommand{\hzmodmap}[2]{\hzmod(#1,#2)}

\newcommand{\hztumod}{\modules_{H\ft[u]}}
\newcommand{\hztmod}{\modules_{H\ft}}
\DeclareMathOperator{\modules}{\bf Mod}
\newcommand{\zmodf}{\zmod^{T}}
\newcommand{\hzmodf}{\hzmod^{T}}

\DeclareMathOperator{\Ab}{\bf Ab}
\DeclareMathOperator{\Vect}{\bf Vect}
\newcommand{\frab}{\Ab_{\evenstar}^{\text{fr}}}

\newcommand{\abstar}{\Ab_{\ast}}
\newcommand{\grvect}{\Vect_{\ast}}

\newcommand{\mot}{\mathbf{M}}
\newcommand{\motZ}{\mathbf{C}}
\newcommand{\vmot}{\mathbf{V}}
\newcommand{\motGen}{\mathbf{G}}

\newcommand{\zlp}[1]{\Z_{(p)}}

\newcommand{\evenstar}{\bullet}

\DeclareMathOperator{\res}{res}

\newcommand{\brk}[1]{\rwHringIdeal^{#1}}
\DeclareMathOperator{\image}{image}

\newcommand{\einfty}{E_{\infty}}
\newcommand{\cat}{\mathcal}

\newcommand{\kk}{2k}

\DeclareMathOperator{\id}{Id}

\makeatletter
\makeatletter
\newcommand{\oset}[3][0ex]{%
  \mathrel{\mathop{#3}\limits^{
    \vbox to#1{\kern-2\ex@
    \hbox{$\scriptstyle#2$}\vss}}}}
\makeatother

\DeclareMathOperator{\alg}{Alg}

\newcommand{\sets}{\mathfrak{Sets}}

\newcommand{\unit}{\mathbf{1}}

\newcommand{\freeRmodule}[2]{\mathscr{F}^{\text{mod}}_{#1}{(#2)}}

\newcommand{\addGrp}{\mathbb{A}}

\newcommand{\tfunc}{\mathbf{T}}

\DeclareMathOperator{\Primitives}{Prim}
\newcommand{\prim}[1]{\Primitives_{(1)}}

\newcommand{\mutn}[1]{MU^{(\!(#1)\!)}}
\newcommand{\versch}{\mathbf{v}}
\newcommand{\partitions}[1]{p(#1)}

\begin{document}

\maketitle

\tableofcontents
\section*{Introduction}

In~\cite{WSWThesis} Steve Wilson established the remarkable fact that
the even spaces 
\[
\umu^{2n}=\Omega^{\infty}\Sigma^{2n}MU
\]
in the complex cobordism spectrum have cell decompositions with only
even dimensional cells.  One consequence of this is that the unstable
$MU$ Adams resolution for a CW complex $X$ with only even cells takes
the form
\[
X \to I_{0}\to I_{1}\to\cdots
\]
in which $I_{n}$ is the zeroth space of a spectrum of the form
\[
\bigvee S^{2n_{\alpha}}\wedge MU.
\]
This leads to an unexpected resolution of $X$ in terms of integer
Eilenberg-MacLane spaces.  

In joint work with Aravind Asok and Jean Fasel, the second author has
explored the hypothesis that the analogue of this holds in motivic
homotopy theory.  This ``Wilson space hypothesis'' leads to
resolutions of certain $\A^{1}$-homotopy types into motivic complexes.
These (conjectural) resolutions have interesting consequences in
algebraic geometry.

Once placed in motivic homotopy theory the question of generality of
the ground ring arises.  If the Wilson Space Hypothesis holds over the
real numbers then the topological realization would be a refinement of
Wilson's theorem to the case of {\em real} homotopy theory.  In fact
this is the case, and part of the purpose of this series of papers is
to establish this refinement.  The methods employed actually lead to
an analogous result for the spectra $\mutn{C_{2^{n}}}$ introduced
in~\cite{MR3505179} and in the third paper in this series we formulate
and establish this result.  Taken together these theorems suggest that
the Wilson Space Hypothesis is something very general, and that one
should seek a less computational approach.

In this paper we consider spaces with a $\zt$-action, with examples in
mind being the complex points of a variety defined over $\R$.  In this
spirit we regard $\C^{n}$ as a $\zt$ space under complex conjugation,
and write
\[
S(\C^{n})\subset D(\C^{n})\subset \C^{n}
\]
for the unit spheres and disks.  The natural analogue of a CW complex
consisting only of even dimensional cells is a $\zt$-space which
decomposes into cells of the form $D(\C^{n})$
\[
X=\coprod D(\C^{n}) /\sim
\]
(with the usual closure finite and weak topology assumptions).  We
will say that a $\zt$-space is {\em $\rho$-cellular} if it is equivariantly
weakly equivalent to a space admitting such a decomposition.  The
Schubert cell decompositions show that the complex points of the real
Grassmannians are $\rho$ cellular, and that the classifying space $\bur$ is $\rho$ cellular.
There is an evident analogue of being {\em $\rho$-cellular} in the category of
$\zt$-spectra, and from the Schubert cell decomposition of
Grassmannians one can conclude that the spectrum $\mur$ of {\em real
bordism} is $\rho$ cellular.

For a real representation $V$ of $\zt$ let $S^{V}$ be its one point
compactification.  Following~\cite{MR3505179} the sign representation
of $\zt$ will be denoted by $\sigma$ and the real regular
representation by $\rho=1+\sigma$.  The representation $\rho$ is
realized as the complex numbers with $\zt$-action given by complex
conjugation, and $S^{\rho}$ with $\cp^{1}$.  This is the reason for
the term {\em $\rho$-cellular}.

Given real representations $V$ and $W$ we denote by $S^{V-W}$ the
spectrum constructed as the smash product of $S^{V}$ with the
Spanier-Whitehead dual of $S^{W}$.  The correspondence $V-W\mapsto
S^{V-W}$ sends sums of representations to smash products of
$\zt$-spectra.  With this convention the symbol $S^{V}$ is being used
for both a space and its suspension spectrum.   Which meaning is
intended will be clear from context.  

The real analogue of Wilson's result is given by the following

\begin{thm}
\label{thm:1}
For every $n\in\Z$ the $\zt$-space 
\[
\umur^{n\rho}=\Omega^{\infty}S^{n\rho}\wedge\mur
\]
is $\rho$-cellular.
\end{thm}

Theorem~\ref{thm:1} is proved in~\cite{hill:_real_wilson_ii}, and is
a consequence of the following apparently weaker result.

\begin{thm}
\label{thm:2}
For any $n\in\Z$ there is a weak equivalence 
\[
H\Zm\wedge \bigvee S^{n_{\alpha}\rho} \approx H\Zm\wedge \umur^{n\rho}.
\]
\end{thm}

The implication that Theorem~\ref{thm:2} implies Theorem~\ref{thm:1}
is not very direct.  Theorem~\ref{thm:2} actually yields the
consequences about unstable Adams resolutions one would wish to derive
from knowing that $\umur^{n\rho}$ is $\rho$-cellular, so it's not clear
which of these two results is most natural to seek in general.

Most of the work of this paper involves setting up a useful
description of the map of Theorem~\ref{thm:2}, with the more refined
statement appearing as Theorem~\ref{thm:26}.  This is carried out
using the language of {\em Hopf rings}.  Hopf rings were introduced by
Ravenel and Wilson in~\cite{RW:HR}, where the homology of the spaces
constituting the complex cobordism spectrum where characterized by a
universal property.  For equivariant and motivic homotopy theory a
refinement of this universal property is needed and
in~\cite{hill18:_univer_hopf} the theory of Hopf rings is developed
and this refined universal property is established.  A summary of the
necessary results is given in \S\ref{sec:universal-hopf-rings}.

Here is a more detailed summary of the contents of this paper.  The
first section simply collects the notation and conventions used
throughout the paper.  In \S\ref{sec:universal-hopf-rings} we recall
the Ravenel-Wilson theory of ``Hopf rings'' and describe our
refinement of the Ravenel-Wilson universal property.
Section\ref{sec:some-algebra-i} contains some mildly technical results
about the Hopf rings of interest in this paper.  In
\S\ref{sec:equiv-r-modul} we turn to $\zt$ equivariant homotopy
theory, recall the theory of {\em real oriented} spectra, and state in
precise terms our equivariant refinement of the result of
Ravenel-Wilson (Theorem~\ref{thm:26}).  The proof of
Theorem~\ref{thm:26} reduces to a question about the mod $2$ homology
of certain fixed point spaces $X^{2k}$ appearing in the real bordism
spectrum.  The basic facts about the $X^{2k}$ are established in
\S\ref{sec:hopf-ring-xbullet}.  Finally, everything is put together in
\S\ref{sec:proof-main-theorem} where the proof of Theorem~\ref{thm:26}
is completed.  The first part of \S\ref{sec:proof-main-theorem}
recollects the entirety of the argument, so, after becoming acquainted
with the basic setup of this paper, the reader may wish to look there
for a further overview.

We prove results about the homology of the spaces $X^{2k}$ by
following the inductive argument of Chan~\cite{MR675011}.    The
homology of these spaces has been computed by Kitchloo and Wilson as
part of a much more general result~\cite[Theorem~1.5]{MR2317067}, and
our main theorem about them (Proposition~\ref{thm:35}) can also 
be deduced from the Kitchloo-Wilson theorem.  

\subsubsection*{Acknowledgements}
\label{sec:acknowledgements}

The authors would like to thank John Francis for the very helpful
suggestion of working in the homotopy category of $H\Zm$-modules.
They would also like to thank Danny Shi and Dylan Wilson for their
interest in this work and for insisting on hearing all of the details.

During the preparation of this manuscript the first author was
supported by NSF grant DMS-1509652 and the second author by 
NSF grant DMS-1510417.  

\numberwithin{equation}{section}
\numberwithin{figure}{section}

\section{Universal Hopf Rings}
\label{sec:universal-hopf-rings}

In this section we recall the Ravenel-Wilson theory of Hopf rings and
describe the universal property characterizing the Hopf ring for
complex cobordism.  Our main purpose is to state a generalization of
this universal property (Theorem~\ref{thm:6}).  This generalization is
proved in~\cite{hill18:_univer_hopf}, and used here to construct the
map~\eqref{eq:20} appearing in the statement of the main theorem of
this paper (Theorem~\ref{thm:26}).

\subsection{Hopf rings and the Ravenel-Wilson theorem}
\label{sec:hopf-rings-ravenel}

In this section we define the notion of a {\em Hopf ring} in a
symmetric monoidal category.   More complete details appear in~\cite{hill18:_univer_hopf}.

\subsubsection{Hopf Rings}
\label{sec:hopf-rings-1}
Suppose that $\cat C=(\cat C,\otimes)$ is a symmetric monoidal
category.  Let $\coalg \cat C$ be the category of counital,
cocommutative, coassociative coalgebras in $\cat C$.  The category $\coalg\cat C$ has
finite products, given on underlying objects in $\cat C$ by the
monoidal product $\otimes$.

\begin{defin}
\label{def:1}
A {\em Hopf ring} in $\cat C$ is a commutative ring object in $\coalg
\cat C$.
\end{defin}

Thus a Hopf ring consists of an object $\hringh\in \coalg\cat C$ equipped
with addition and multiplication laws 
\begin{equation}
\label{eq:1}
\begin{aligned}
+:\hringh\times \hringh &\to \hringh \\
\times:\hringh\times \hringh &\to \hringh \\
\end{aligned}
\end{equation}
and additive and multiplicative units
\begin{equation}
\label{eq:2}
[0], [1]:\unit\to \hringh
\end{equation}
satisfying the axioms of a commutative ring.   Note that the products
in~\eqref{eq:1} are formed in $\coalg \cat C$, so correspond to the
monoidal product $\otimes$ in $\cat C$, and the maps~\eqref{eq:2}
and~\eqref{eq:1} are coalgebra maps.

By the Yoneda lemma, giving a Hopf ring structure on $\hringh$ is equivalent
to lifting the functor 
\[
\hringh(\slot)=\coalg\cat
C(\slot,\hringh)
\]
to one taking values in commutative rings (see \cite{hill18:_univer_hopf}).

\begin{eg}
\label{eg:1} The most basic example is when $\cat C=\sets$, with the
symmetric monoidal structure given by the Cartesian product.  Since
every set has a unique counital coalgebra structure, the forgetful
functor $\coalg\cat C\to \sets$ is an equivalence of categories.  A
Hopf ring in $\sets$ is just a commutative ring.
\end{eg}

\begin{eg}
\label{eg:2} Similarly, a Hopf ring in the homotopy category
$\ho\spaces$ of spaces and cartesian products is a ring object in
$\ho\spaces$.  For instance if $E$ is a homotopy commutative ring
spectrum then the space $\underline{E}^{0}=\Omega^{\infty}E$ inherits
the structure of a Hopf ring in $\ho\spaces$.
\end{eg}

\begin{eg}
\label{eg:3} 
Suppose that $k$ is a commutative ring and let
$\gmodules{k}$ be the symmetric monoidal category of graded left
$k$-modules, and tensor products with the symmetry given by the usual
Koszul sign convention.  In concrete terms, a Hopf ring $\hringh$ over
$\gmodules{k}$ is a graded $k$-module equipped with maps
\begin{align*}
\psi &:\hringh \to \hringh\otimes \hringh \\
\epsilon &:\hringh \to k,
\end{align*}
making it into a coalgebra, coalgebra maps
\begin{align*}
\ast &:\hringh\otimes \hringh\to \hringh \\
\circ &:\hringh\otimes \hringh\to \hringh
\end{align*}
corresponding to addition and multiplication, and
constants~\eqref{eq:2}.  The constants may be identified with elements
\begin{align*}
[0] & \in \hringh \\
[1] &\in \hringh.
\end{align*}
The product $\ast$ makes each $\hringh$ into a commutative, cocommutative
Hopf algebra over $k$.  To keep the notation simple we will denote it
by ordinary multiplication.  The multiplicative unit for this ring
structure is the element $[0]$, so one has
\[
1=[0]\in \hringh.
\]   
\end{eg}

If $F:\cat C\to \cat D$ is a symmetric monoidal functor and $\hringh$ is a
Hopf ring in $\cat C$ then $F(\hringh)$ is a Hopf ring in $\cat D$.

\begin{eg}
\label{eg:4}
Suppose that $R$ is an $E_{\infty}$ ring spectrum and let $\rmod{R}$
be the homotopy category of left $R$-modules, regarded as an additive
category.   Equipped with the monoidal product 
\[
M\underset{R}{\wedge}N
\]
the category $\rmod{R}$ becomes a symmetric monoidal category.    The
functor  
\begin{gather*}
\freeRmodule{R}{X} = R\wedge X_{+}
\end{gather*}
is a symmetric monoidal functor from $\ho\spaces$ to $\rmod{R}$.  If
$E$ is a homotopy commutative ring spectrum then
$\freeRmodule{R}{\Omega^{\infty} E}=R\wedge \Omega^{\infty}E_{+}$
is a Hopf ring in $\rmod{R}$.
\end{eg}

\begin{eg}
\label{eg:5} If $k$ is a ring, $E$ is a homotopy commutative ring
spectrum, and $\underline{E}^{0}=\Omega^{\infty}E$ has the property
that
\[
H_{\ast}(\underline{E}^{0};k)
\]
is flat over $k$ then $H_{\ast}(\underline{E}^{0};k)$ is a Hopf ring in
$\gmodules{k}$.
\end{eg}

If $F:\cat C\to \cat D$ is a symmetric monoidal functor, $\hringh$ is a
Hopf ring in $\cat C$, and $X\in\coalg\cat C$ is a coalgebra, then $F$
provides a ring homomorphism 
\[
\hringh(X) \to (F\hringh)(FX)
\]
where, as in the discussion just before Example~\ref{eg:1}, $\hringh(X)$
denotes the functor $\coalg\cat C(X,\hringh)$ and $(F\hringh)(FX)= \coalg\cat D(F\hringh,FX)$.

\begin{eg}
\label{eg:6} 
In the situation of Example~\ref{eg:4}, if $X$ is a
space and $E$ is a homotopy commutative ring spectrum then the functor
$\freeRmodule{R}{\slot}$ gives a ring homomorphism
\[
\underline{E}^{0}(X) = E^{0}(X) \to
\freeRmodule{R}{\underline{E}^{0}}(\freeRmodule{R}{X}) =
\coalg\rmod{R}(R\wedge X_{+}, R\wedge \underline{E}^{0}_{+}),
\]
which may be interpreted as a natural transformation of ring valued
functors of $X$.  
\end{eg}

\subsubsection{Hopf rings over a ring}
\label{sec:hopf-rings-over-1}

Suppose that $A$ is an ordinary commutative ring.

\begin{defin}
\label{def:2}
A {\em Hopf ring over $A$} in $\cat C$ is is a commutative
$A$-algebra in $\coalg\cat C$.
\end{defin}

Equivalently, a Hopf ring over $A$ in $\cat C$ a representable functor
from $\coalg\cat C$ to the category $\alg_{A}$ of commutative
$A$-algebras.  In addition to the ring structure, a Hopf ring $\hringh$ over
$A$ has constants
\[
[a]:\unit\to \hringh \qquad a\in A
\]
which, together with the
addition, multiplication and units satisfy the axioms of an
$A$-algebra.

If $F:\cat C\to \cat D$ is symmetric monoidal, and $\hringh$ is a Hopf ring
over $A$ in $\cat C$ then $F\hringh$ is a Hopf ring over $A$ in $\cat D$.

\begin{eg}
\label{eg:7} If $E$ is a ring spectrum then then
$\underline{E}^{0}=\Omega^{\infty}E$ is a Hopf ring over the ring
$E^{0}=E^{0}(\text{pt})$ in the homotopy category of spaces.  By
functoriality, the free $R$-module
$\freeRmodule{R}{\underline{E}^{0}}$ is a Hopf ring over $E^{0}$ in
$\rmod{R}$.
\end{eg}

\subsubsection{Evenly graded Hopf rings}
\label{sec:evenly-grade-hopf}

For a category $\cat D$ let $\cat D^{\bullet}$ be the category of
evenly graded objects of $\cat D$.  Thus an object $X^{\bullet}\in
\cat D^{\bullet}$ consists of objects $X^{m}\in\cat D$ for every {\em even}
integer $m\in\Z$.  When $\cat C=(\cat C,\otimes)$ is
a symmetric monoidal category we will write
\[
\coalg \cat C^{\bullet} = \big(\coalg \cat C)^{\bullet}
\]
for the category of evenly graded objects in $\coalg \cat C$.  We will
also employ the notation $X^{\bullet}$ for an evenly graded object of
$\cat D$, and use the convention $X^{m}=X_{-m}$.     

Suppose that $\rbullet{A}$ is an evenly graded commutative ring.

\begin{defin}
\label{def:3}
A {\em Hopf ring over $\rbullet{A}$} in $\cat C$ is an evenly graded
commutative $\rbullet{A}$-algebra in $\coalg\cat C^{\bullet}$.  
\end{defin}

Thus a Hopf ring over $\rbullet{A}$ in $\cat C$ is an evenly graded
object $\hringh^{\bullet}\in \coalg\cat C^{\bullet}$, equipped with addition
and multiplication maps 
\begin{align*}
+:\hringh^{2n}\times \hringh^{2n}&\to \hringh^{2n} \\
\times : \hringh^{n}\times \hringh^{2m}&\to \hringh^{2n+2m},
\end{align*}
constants 
\[
[a] :\unit\to \hringh^{2n}
\]
for every $a\in A^{2n}$,  including the additive units 
\[
[0]:\unit\to \hringh^{2n},
\]
satisfying the axioms of a graded commutative $\rbullet{A}$-algebra.
We have restricted to {\em evenly} graded objects in order that there
be no confusion with the Koszul sign convention.  On evenly graded
objects, commutativity is just ordinary commutativity.
See~\cite[\S1]{RW:HR} and~\cite{hill18:_univer_hopf} for more
explicated details.

\begin{eg}
\label{eg:8} 
For a spectrum $E$ set
\[
\underline{E}^{m} = \Omega^{\infty}\Sigma^{m}E,
\]
and, ignoring the odd values of $m$, consider the evenly graded space
$\underline{E}^{\bullet}$.  If $E$ is a homotopy commutative ring
spectrum then $\underline{E}^{\bullet}$ is a Hopf ring in spaces over
the evenly graded ring $\rbullet{E}$, in which
$E^{2n}=E^{2n}(\text{pt})$.
\end{eg}

\begin{eg}
\label{eg:9} By functoriality, if $R$ is an $E_{\infty}$ ring then
$\freeRmodule{R}{\underline{E}^{\bullet}}=R\wedge
\underline{E}^{\bullet}_{+}$ is an evenly graded Hopf ring over
$E^{\bullet}$ in $\rmod{R}$.
\end{eg}

\begin{eg}
\label{eg:10}
If $k$ is a ring, and for each $n\in\Z$, the graded
$k$-module 
\[
H_{\ast}(\underline{E}^{2n};k)
\]
is flat, then $H_{\ast}(\underline{E}^{\bullet};k)$ is a Hopf ring
over $E^{\bullet}$ in $\gmodules{k}$.   
\end{eg}

\begin{notation}
\label{def:4} The category $\hopfrings_{\rbullet{R}}(\cat C)$ is the
category of evenly graded Hopf rings over $\rbullet{R}$ in $\cat C$,
and Hopf ring homomorphisms over $\rbullet{R}$.
\end{notation}

\subsection{The Ravenel-Wilson theory}
\label{sec:ravenel-wilson-hopf-1}

\subsubsection{Examples from complex oriented cohomology theories}
\label{sec:exampl-from-compl}

We are especially interested in the evenly graded Hopf rings over
$\rbullet{MU}$ with $MU^{2n}=MU^{2n}(\text{pt})$, and the additional
structures present in those that arise from a complex oriented
cohomology theory.  Suppose that $(E,x)$ is a complex oriented
homotopy commutative ring spectrum and let
\begin{equation}
\label{eq:3}
MU^{\bullet}\to \rbullet{E}
\end{equation}
be the map classifying the associated formal group law.   As in
\S\ref{sec:evenly-grade-hopf}, the evenly graded space
$\underline{E}^{\bullet}$ is an evenly graded Hopf ring over
$E^{\bullet}$.   Restricting scalars along~\eqref{eq:3} makes
$\underline{E}^{\bullet}$ into a Hopf ring over $\rbullet{MU}$.
If $R$ is an $E_{\infty}$ ring then, by Example ~\ref{eg:9}, the $R$-module
\[
\freeRmodule{R}{\underline{E}^{\bullet}}
\]
is a Hopf ring over $\rbullet{MU}$ in the homotopy category
$\rmod{R}$ of left $R$-modules.

By definition of the formal group law, the complex orientation
$x\in E^{2}(\cp^{\infty})$ satisfies 
\begin{equation}
\label{eq:4}
\begin{aligned}
\beta_{0}^{\ast}x &= 0  \\
\mu^{\ast}x &= \sum a_{ij}x^{i}y^{j} ,
\end{aligned}
\end{equation}
where 
\begin{equation}
\label{eq:62}
\mu:\cp^{\infty}\times \cp^{\infty}\to \cp^{\infty}
\end{equation}
is the map classifying the tensor product of the two universal line
bundles, the $a_{ij}\in E^{-2(i+j)}(\text{pt})$ are the coefficients
of the formal group law, and
\[
\beta_{0}:\{\text{pt}\}\to\cp^{\infty}
\]
is the homotopy class of the inclusion of a point.  The
identity~\eqref{eq:4} may be interpreted in the degree $2$ component
of the ring $\underline{E}^{\bullet}(\cp^{\infty}\times \cp^{\infty})$
associated to the Hopf ring $\underline{E}^{\bullet}$ over
$\rbullet{MU}$ in $\ho\spaces$.  By functoriality (Example~\ref{eg:6})
the coalgebra map
\[
\tilde x= \freeRmodule{R}{x} : \freeRmodule{R}{\cp^{\infty}}\to
\freeRmodule{R}{\underline{E}^{2}} 
\]
satisfies 
\begin{equation}
\label{eq:5}
\begin{aligned}
\beta_{0}^{\ast}(\tilde x) &= [0] \\
\mu^{\ast}\tilde x &= \sum a_{ij} \tilde x^{i} \tilde y^{j}
\end{aligned}
\end{equation}
in which we have written $\beta_{0}=\freeRmodule{R}{\beta_{0}}$ and
$\mu = \freeRmodule{R}{\mu}$.  To further simplify the notation we
will drop the tilde and write $x=\freeRmodule{R}{x}$, etc in the above
identity.

If $k$ is a ring and for each $n\in\Z$,
$H_{\ast}(\underline{E}^{2n};k)$ is a flat $k$-module, then
$H_{\ast}(\underline{E}^{\bullet};k)$ is a Hopf ring in $\gmodules{k}$
over $\rbullet{MU}$, equipped with an element
\[
x:H_{\ast}(\cp^{\infty};k)\to H_{\ast}(\underline{E}^{2};k)
\]
satisfying~\eqref{eq:5}.

\subsubsection{The Ravenel-Wilson Hopf ring}
\label{sec:ravenel-wilson-hopf-2}

In~\cite{RW:HR} Ravenel and Wilson construct a universal Hopf ring
over $\rbullet{MU}$ equipped with an element $x$ satisfying the
relation~\eqref{eq:5}.  

To make this precise, let $\abstar=\gmodules{\Z}$ be the symmetric
monoidal category of graded abelian groups (with the Koszul sign
convention).  Denote by $\hopfRingsMU$ the category of Hopf rings over
$\rbullet{MU}$ in $\abstar$.  For brevity, we will refer to the
objects of $\hopfRingsMU$ as {\em evenly graded Hopf rings over
$\rbullet{MU}$}, or sometimes as just Hopf rings.

Let $\addGrp^{1}$ be the abelian group object in $\coalg\abstar$ given
by the Hopf algebra $H_{\ast}(\cp^{\infty};\Z)$.  Thus $\addGrp^{1}$
is the coalgebra
\[
\Z\{\beta_{0},\beta_{1}, \cdots \}\qquad |\beta_{i}|=2i,
\]
with coproduct given by 
\[
\beta_{n}\mapsto \sum_{i+j=n} \beta_{i}\otimes \beta_{j}.
\]
The product is induced by the map~\eqref{eq:62}, and
given explicitly by
\[
\beta_{i}\beta_{j}=\binom{i+j}{i} \beta_{i+j}.   
\]
For $n\ge 0$ set
\[
\A^{n}=(\addGrp^{1})^{n} = H_{\ast}\big((\cp^{\infty})^{n};\Z\big).
\]
Write 
\[
\pi_{i}:\A^{n}\to \addGrp^{1}
\]
for projection to the $i^{\text{th}}$ factor, and 
\[
\mu:\A^{2}\to \addGrp^{1}
\]
for the group law described above.  Since our coalgebras are all
counital, the coalgebra map
\[
\beta_{0}:\A^{0}\to \A^{1}
\]
sending $1\in\Z$ to $\beta_{0}$ is the unique  coalgebra
map from $\A^{0}$ to $\A^{1}$.

\begin{defin}
\label{def:5} 
Suppose that $\hringh^{\bullet}$ is an evenly graded
Hopf ring over $\rbullet{MU}$ in $\coalg\abstar$.  An {\em additive
curve} of $\hringh^{\bullet}$ is an element
\[
x\in \hringh^{2}(\A^{1})
\]
satisfying 
\begin{align*}
\beta_{0}^{\ast}x &= [0]\in \hringh^{2}(\unit) \\
\mu^{\ast}x &= \sum a_{ij}x^{i}y^{j}.
\end{align*}
\end{defin}

\begin{defin}
\label{def:6} 
The {\em Ravenel-Wilson Hopf ring} is the evenly graded Hopf ring
$\hringmurw^{\bullet}$ over $\rbullet{MU}$ characterized by the
following universal property: to give a map
\[
\hringmurw^{\bullet}\to \hringh^{\bullet}
\]
in $\hopfRingsMU$ is equivalent to giving an element $x\in
\hringh^{2}(\addGrp^{1})$ having the properties 
\begin{equation}
\label{eq:6}
\begin{aligned}
\beta_{0}^{\ast}x &= [0]\in \hringh^{2}(\Z) \\
\mu^{\ast}x &= \sum a_{ij}x^{i}y^{j}
\end{aligned}
\end{equation}
in which the $a_{ij}\in MU^{-2(i+j)}$ are the coefficients of the
universal formal group law.
\end{defin}

Put differently, the Hopf ring $\hringmurw^{\bullet}$ corepresents
the functor on $\hopfRingsMU$ sending $\hringh^{\bullet}$  to the set of
additive curves of $\hringh^{\bullet}$.  

A proof of the existence of $\hringmurw^{\bullet}$ is sketched
in~\cite{RW:HR}, and given in detail
in~\cite{MR1631257},~\cite{MR2559638} and~\cite{hill18:_univer_hopf}.
As for the structure of $\hringmurw^{\bullet}$, Ravenel and Wilson showed

\begin{thm}
\label{thm:3} 
For each $n$ the algebra $\hringmurw^{2n}$ is a
polynomial algebra over the group ring $\Z[MU^{2n}]$ on generators of
positive even degree.  In particular, $\hringmurw^n$ is free when
regarded as a graded abelian group.
\end{thm}

As mentioned in the introduction, the following theorem of
Wilson started this whole investigation.

\begin{thm}[\cite{WSWThesis}] 
\label{thm:4} 
For every integer $n$, the graded abelian group
$H_{\ast}(\umu^{n};\Z)$ is free.  It is the zero group when $n$ is
even and $\ast$ is odd.
\end{thm}

By the discussion at the end of \S\ref{sec:exampl-from-compl},
Theorem~\ref{thm:4} implies that the graded abelian groups
$H_{\ast}(\umu^{\bullet};\Z)$ form a Hopf ring over $\rbullet{MU}$ in
$\abstar$, equipped with an additive curve $x\in H_{\ast}(\umu^{2})(\A^{1})$.

\begin{thm}[\cite{RW:HR}, Corollary~4.7] 
\label{thm:5} 
The map of Hopf rings over $\rbullet{MU}$
\[
\hringmurw^{\bullet} \to H_{\ast}(\umu^{\bullet})
\]
classifying the additive curve $x$ is an isomorphism.
\end{thm}

\begin{rem}
\label{rem:1} 
In fact, in~\cite{WSWThesis} Wilson also determined the
ring structure of $H_{\ast}\underline{BP}^{n}$ and hence
$H_{\ast}\umu^{n}$ for all $n$.  One purpose of~\cite{RW:HR} was to
give new proof of this result organized by the language of Hopf rings.
\end{rem}

In~\cite{RW:HR} the proofs of Theorems~\ref{thm:3},~\ref{thm:4} and~\ref{thm:5}
are interleaved, and in fact the result is first proved for homology
with coefficients in a field and then the general result is
deduced.  Without either working over a field, or first assuming
Wilson's Theorem~\ref{thm:4}, one does not know in advance that the abelian
groups $H_{\ast}\umu^{\bullet}$ are even coalgebras, much less Hopf
rings.  In~\cite{hill18:_univer_hopf} the existence of
$\hringmurw^{\bullet}$ is established purely algebraically, as is
Theorem~\ref{thm:3}.  Using Theorem~\ref{thm:3} one can establish a
stronger universal property of the Ravenel-Wilson Hopf ring and use it
to make the desired map above.  This is the subject of the next
section.

\subsubsection{Generalizing the Ravenel-Wilson universal property}
\label{sec:gener-raven-wils}

Let $\frab$ be the symmetric monoidal category of evenly graded free
abelian groups and tensor products.  By Theorem~\ref{thm:3} the Hopf
ring $\hringmurw^{\bullet}$ can be regarded as a Hopf ring over
$\rbullet{MU}$ in $\frab$.  Suppose that $(\cat C,\otimes)$ is a
strongly additive symmetric monoidal  category in the sense of
Appendix~\S\ref{sec:notat-basic-noti}, and write $\hopfRingsMU(\cat C)$ for the
category of Hopf rings over $\rbullet{MU}$ in $\cat C$.  If
\[
\motGen:\frab\to \cat C
\]
is an additive symmetric monoidal functor then
$\motGen(\hringmurw^{\bullet})$ is a Hopf ring over $\rbullet{MU}$ in
$\cat C$.   Using $\motGen$ we make the following definition.

\begin{defin}
\label{def:7} 
Suppose that $\hringh^{\bullet}$ is an evenly graded
Hopf ring over $\rbullet{MU}$ in $\cat C$.  An {\em additive curve} in
$\hringh^{\bullet}$ is a coalgebra map
\[
x:\motGen(\A^{1})\to \hringh^{2}
\]
satisfying 
\begin{align*}
\motGen(\beta_{0})^{\ast}x &= [0] \\
\motGen(\mu)^{\ast}x &= \sum a_{ij}x^{i}y^{j}.
\end{align*}
\end{defin}

Let
\[
\tfunc:\hopfRingsMU(\cat C) \to \sets
\]
be the functor sending $\hringh^{\bullet}$ to the set of additive
curves in $\hringh^{\bullet}$.  Since $\motGen$ is additive,
$\motGen(x_{\text{univ}})$ is an element of
$\tfunc(\motGen(\hringmurw^{\bullet}))$.  The main result
of~\cite{hill18:_univer_hopf} is

\begin{thm}
\label{thm:6}
With the above notation, the functor 
\[
\tfunc:\hopfRingsMU(\cat C)\to \sets
\]
is represented by $(\motGen(\hringmurw^{\bullet}),
\motGen(x_{\text{univ}}))$.   
\end{thm}

Phrased more succinctly, the image of the Ravenel-Wilson Hopf ring
under a strongly additive symmetric monoidal functor defined on free abelian
groups enjoys the same universal property as the Ravenel-Wilson Hopf
ring.  

Here is a simple application of Theorem~\ref{thm:6}.  There is a
unique strongly additive symmetric monoidal functor
\begin{equation}
\label{eq:7}
\motZ:\frab\to \zmod
\end{equation}
satisfying 
\[
\motZ(\Z[2n])= H\Z\wedge S^{2n}.
\]
In fact $\motZ$ is an equivalence between $\frab$ and the smallest
full subcategory of $\zmodf\subset \zmod$ containing the $S^{2n}\wedge H\Z$ and
closed under arbitrary coproducts.  If $(E,x)$ is a homotopy
commutative complex oriented cohomology theory then by
\S\ref{sec:exampl-from-compl}, the evenly graded coalgebra
$\freeRmodule{H\Z}{\underline{E}^{\bullet}}= H\Z\wedge
(\underline{E}^{\bullet})_{+}$ becomes a Hopf ring over $\rbullet{MU}$
in $\coalg\zmod$ equipped with an additive curve.  By
Theorem~\ref{thm:6} there is a unique Hopf ring map
\[
\motZ(\hringmurw) \to \freeRmodule{H\Z}{\underline{E}^{\bullet}}
\]
classifying the additive curve.  When $E=MU$  this can be used to
define the map of Theorem~\ref{thm:5} directly.

\begin{rem}
\label{rem:3}
The inverse of the equivalence $\frab\to \zmodf$ sends $M$ to
$\pi_{\bullet}M$.   By definition then,  one has
\[
\motZ(\A^{k})=\freeRmodule{H\Z}{(\cp^{\infty})^{k}}=
H\Z\wedge(\cp^{\infty})^{k}_{+}\mathrlap{\ .}
\]
\end{rem}

While this application of Theorem~\ref{thm:6} is relatively minor, in
equivariant and motivic homotopy theory it helps significantly.

\section{Verschiebung and decomposition}
\label{sec:some-algebra-i}

Our aim in this section is to set up a certain tensor product
decomposition (Proposition~\ref{thm:12}) that will aid the proof of
our main result (Theorem~\ref{thm:26}).  We also recall an important
technical result of Ravenel and Wilson (Proposition~\ref{thm:13}).

\subsection{Connected components}
\label{sec:connected-components}

Suppose that $\cat C=\gmodules{k}$ is the category of graded modules
over a commutative ring $k$, made into a symmetric monoidal category
with the tensor product and the Koszul sign rule.  

\begin{defin}
\label{def:8} 
A coalgebra $C\in \coalg\cat C$ is {\em
$(-1)$-connected} if $C_{i}=0$ for $i<0$.  A coalgebra is {\em
connected} if it is $(-1)$-connected and if the counit map
$\epsilon:C\to k$ is an isomorphism.
\end{defin}

When $C$ is $(-1)$ connected the projection map $p:C\to C_{0}$ is a
coalgebra map.

\begin{defin}
\label{def:9}
A {\em grouplike element of $C$} is a coalgebra map $\unit\to C$.
\end{defin}

\begin{rem}
\label{rem:4}
Since our maps of coalgebras are counital, the composition $\unit\to
C\xrightarrow{\epsilon}{} \unit$ is the identity map.
\end{rem}

\begin{defin}
\label{def:17}
The set $\pi_{0}C$ is the set
\[
C(\unit)=\coalg\cat C(\unit,C)
\]
of grouplike elements of $C$.
\end{defin}

A grouplike element $\unit\to C$ is uniquely determined by the image
of $1\in k$ which can be any element $x\in C_{0}$ with $\epsilon(x)=1$ and
$\psi(x)=x\otimes x$ where $\psi:C\to C\otimes C$ is the coalgebra
structure map.   

\begin{eg}
\label{eg:12} 
The functor sending a set $T$ to the free $k$-module
$k\{T \}$ is symmetric monoidal, giving $k\{T \}$ the structure of a
coalgebra.  When $k$ has no idempotents other than $0$ and $1$, the
inclusion of the generators $T\to k\{T \}$ is an isomorphism of $T$
with the set of grouplike elements of $k\{T \}$.
\end{eg}

There is a canonical coalgebra map 
\begin{equation}
\label{eq:8}
k\{\pi_{0}C\}=\bigoplus_{\pi_{0}C} k\to C_{0}.
\end{equation}

\begin{defin}
\label{def:10} A coalgebra $C\in \gmodules{k}$ is {\em spacelike} if
$C$ is $(-1)$-connected, and the
map~\eqref{eq:8} is an isomorphism.
\end{defin}

\begin{eg}
\label{eg:16} 
Let $k$ be a commutative ring with no idempotents other
than $0$ and $1$.  If $X$ is a space and $H_{\ast}(X;k)$ is flat over
$k$, then $H_{\ast}(X;k)$ is a coalgebra in $\gmodules{k}$.  The
coalgebra map induced by the inclusion of a point $x\in X$ defines a
grouplike element of $H_{\ast}(X;k)$, depending only on the path
component of $x$.  The resulting map $\pi_{0}X\to
\pi_{0}H_{\ast}(X;k)$ is a bijection and $H_{\ast}(X;k)$ is spacelike.
\end{eg}

Suppose that $C$ is a $(-1)$-connected coalgebra and $b\in\pi_{0}C$ is a grouplike
element.   

\begin{defin}
\label{def:11}
The {\em connected component of $C$}  containing $b$ is the
sub coalgebra $C'=C'_{b}\subset C$ consisting of elements $x\in C$ whose image
under 
\[
C\to C\otimes C\xrightarrow{\,\id\otimes p\,} C\otimes C_{0}
\]
is $x\otimes b$.
\end{defin}

That $C'$ actually is a subcoalgebra follows easily from the
coassociativity law, and the fact that the map $b:\unit\to C$ is split
by the counit and so the inclusion of a $k$-module summand.

\begin{rem}
\label{rem:11}
The coalgebra $C'_{b}$ fits into a pullback diagram 
\begin{equation}
\label{eq:41}
\xymatrix{
C'_{b}  \ar[r]\ar[d]  & C
\ar[d] \\
\unit  \ar[r]        & C_{0}
}
\end{equation}
in $\coalg \gmodules{k}$.   
\end{rem}

\begin{eg}
\label{eg:13} 
Suppose $X$ is a pointed space for which $H_{\ast}(X;k)$
is flat over $k$.  As in Example~\ref{eg:16}, the homology class $b$
of the base point is a grouplike element of $H_{0}X$, and the
coalgebra connected component
\[
H_{\ast}(X;k)' = H_{\ast}(X;k)_{b}'
\]
is the homology of the connected component $X'\subset X$ of $X$
containing the base point.
\end{eg}

We highlight one piece of notation used in the above example.

\begin{notation}
If $X$ is a pointed space, with base point $x\in X$, then $X'\subset
X$ is the component of $X$ containing the base point.
\end{notation}

\subsection{Hopf algebras}
\label{sec:hopf-algebras}

\begin{defin}
\label{def:18}
Let $\cat C$ be a symmetric monoidal category.  A {\em Hopf algebra in
$\cat C$} is an abelian group object $\mathcal A\in \coalg\cat C$.
\end{defin}

When $\cat C=\gmodules{k}$ a Hopf algebra $\mathcal A$ in $\cat C$ can
be described as a coalgebra $\mathcal A$, equipped with a coalgebra
map
\[
\mathcal A\otimes \mathcal A\to \mathcal A
\]
and an additive unit $[0]:\unit \to \mathcal A$ making $\mathcal A$
into a commutative ring.  This data gives
$\coalg\gmodules{k}(\mathcal A,\mathcal A)$ the structure of an abelian
monoid and is subject to the further condition that this abelian
monoid structure forms an abelian group.  This condition is
equivalent to the existence of a (necessarily unique) {\em antipode}
$\chi:A\to A$, which is both a coalgebra and an algebra map, making
the diagram
\[
\xymatrix{
A  \ar[r]\ar[d]_{[0]\circ \epsilon}  & A\otimes A
\ar[d]^{1\otimes \chi} \\
A    & A\otimes A \ar[l]
}
\]
commute.

\begin{notation}
If $\mathcal A$ is a $(-1)$-connected Hopf algebra in $\gmodules{k}$
then $\mathcal A'\subset \mathcal A$  is the connected component of
the additive unit $[0]$.
\end{notation}

\begin{defin}
\label{def:16} 
A Hopf algebra $\mathcal A$ in $\gmodules{k}$ is {\em
spacelike} if the underlying coalgebra is spacelike.
\end{defin}
	
If $\mathcal A$ is a Hopf algebra in $\cat C$ then, by definition, the
set $\pi_{0}\mathcal A=\mathcal A(\unit)$ is an abelian group.

\begin{eg}
\label{eg:15} In the situation of Example~\ref{eg:12}, if $T$ is an
abelian group then the group algebra $k[T]$ is a Hopf algebra in
$\gmodules{k}$, and the isomorphism
\[
T\to \pi_{0}k[T]
\]
is an isomorphism of abelian groups.
\end{eg}

Note that if $\mathcal A$ is a Hopf algebra in $\gmodules{k}$ then the
inclusion $\mathcal A_{0}\to A$ is a map of Hopf algebras.  If
$\mathcal A$ is $(-1)$-connected this map is split by the Hopf algebra
map $\mathcal A\to \mathcal A_{0}$.

\begin{prop}
\label{thm:7} 
If $\mathcal A$ is a $(-1)$-connected Hopf algebra in
$\gmodules{k}$ then the map
\[
\mathcal A_{0}\otimes \mathcal A'\to \mathcal A
\]
is an isomorphism of Hopf algebras.   In particular if 
If $\mathcal A$ is  spacelike then the map 
\[
k[\pi_{0}\mathcal A]\otimes \mathcal A'\to \mathcal A
\]
is an isomorphism.
\end{prop}

\begin{pf}
Since the forgetful functor from  Hopf algebras to coalgebras creates
limits Remark~\ref{rem:11} implies that the square 
\[
\xymatrix{
\mathcal A'_{b}  \ar[r]\ar[d]  & \mathcal A
\ar[d] \\
\unit  \ar[r]        & \mathcal A_{0}
}
\]
is a pullback square of Hopf algebras.   This implies that if $\mathcal
B$ is any Hopf algebra then 
\[
0 \to \hopfalg (\mathcal B,\mathcal A')
\to \hopfalg (\mathcal B,\mathcal A)
\to \hopfalg(\mathcal B,\mathcal A_{0})\to 0
\]
is a split short exact sequence in the category
$\hopfalg=\hopfalg\gmodules{k}$ of abelian group objects in
$\coalg\gmodules{k}$, and so gives an isomorphism
\[
\hopfalg(\mathcal B,\mathcal A)\approx \hopfalg(\mathcal B,\mathcal A')\times
\hopfalg(\mathcal B,\mathcal A_{0}) = \hopfalg(\mathcal B,A'\otimes A_{0}).
\]
The claim follows.
\end{pf}

\subsection{Weight $k$ curves}
\label{sec:weight-k-curves}

Let $\grvect$ be the symmetric monoidal category of $\Z$-graded vector
spaces over $\ft$.  There is a unique coproduct preserving functor
\[
\vmot:\grvect\to \hztmod
\]
sending $\ft[k]$ to $\hft\wedge S^{\kk}$.   The following is straightforward

\begin{prop}
\label{thm:8}
The functor
\begin{align*}
\vmot:\grvect &\to \hztmod
\end{align*}
is an equivalence of strongly additive symmetric monoidal categories, with
inverse given by $\pi_{\ast}$. \qed
\end{prop}

For $k\in \Z$ let $\Phi^{k}:\frab\to \grvect$ be the functor
given by
\[
\Phi^{k}(A)_{m}= \begin{cases}
A_{2n}\otimes \ftwo & \text{if } m = kn \\
0 &\text{otherwise.}
\end{cases}
\]
For each $k$, the functor $\Phi^{k}$ is symmetric monoidal and so
induces functors 
\begin{align*}
\Phi^{k}:\coalg\frab &\to \coalg\grvect \\
\Phi^{k}:\hopfRingsMU(\frab) &\to \hopfRingsMU(\coalg\grvect).
\end{align*}
For $\hringh^{\bullet}\in\hopfRingsMU(\frab)$ and $a\in MU^{2n}$,
the scalar $[a]:\ftwo\to \Phi^{k}\hringh^{2n}$ is obtained by
applying $\Phi^{k}$ to the scalar $[a]:\Z\to \hringh^{2n}$.

For $n\ge 0$ let
\[
\A^{n}(k)=\Phi^{k}(\A^{n})\in \coalg\grvect.
\]
Since $\A^{0}(k)=\A^{0}(0)=\ftwo$ we will just write
$\A^{0}$ for this object.

\begin{defin}
\label{def:12}
An {\em additive curve of weight $k$} in a Hopf ring $\hringh^{\bullet}\in
\hopfRingsMU(\grvect)$  is an element 
\[
x\in \hringh^{2}(\A^{1}(k)).
\]
satisfying 
\begin{align*}
\beta_{0}^{\ast}x &= [0] \\ 
\mu^{\ast}x &= \sum a_{ij}x^{i}y^{j}. 
\end{align*}
\end{defin}

The following is immediate from the definition of $\hringmurw^{\bullet}$.
\begin{prop}
\label{thm:9}
The Hopf ring $\Phi^{k}\hringmurw^{\bullet}$ corepresents the functor
sending a Hopf ring $\hringh^{\bullet}\in \hopfRingsMU(\grvect)$
to the set of additive curves of weight $k$ in
$\hringh^{\bullet}$. \qed
\end{prop}

\subsection{Verschiebung}
\label{sec:verschiebung}

The functor sending a graded vector space $f:V\in\grvect$ to the
graded vector
space $V^{\phi}$ with
\begin{align*}
V^{\phi}_{2n} &=V_{n} \\
V^{\phi}_{2n+1} &=0,
\end{align*}
and a map $f:V\to W$ to the map $f^{\phi}_{2n}=f_{n}$ is symmetric
monoidal.   It therefore induces a product preserving functor
\[
(\slot)^{\phi}:\coalg \grvect\to \coalg \grvect.
\]
The {\em Verschiebung} is a natural transformation
\[
\versch:C\to C^{\phi}
\]
of functors on $\coalg\grvect$.  To define it, note that the map
\[
v \to v\otimes v
\]
gives an isomorphism of $C_{n}$ with the Tate cohomology 
\begin{equation}
\label{eq:9}
\big(\ker (1-\tau)\big)/\big(\image
(1+\tau)),
\end{equation}
where $\tau:(C\otimes C)_{2n}\to (C\otimes C)_{2n}$ is defined by
\[ 
\tau(x\otimes y) = y\otimes x.
\]
For $x\in C_{2n}$, The Verschiebung $\versch(x)$ is the element of $C_{n}$
corresponding to the image of the coproduct $\Delta(x)$
in~\eqref{eq:9}.

While the Verschiebung $\versch$ is always additive in the sense that
\[
\versch(x+y)=\versch(x)+\versch(y),
\]
it is not linear over a general ground field $k$, but rather satisfies
\[
\versch(\lambda^{2}x) = \lambda \versch(x).
\]
In this paper we are working over $\ftwo$, where $\versch$ is in fact
linear.

\begin{rem}
\label{rem:5}
There is a Verschiebung map in any characteristic $p>0$ given by the
modification of the above in which $\tau$ is replace by the cyclic
permutation $C^{\otimes p}\to C^{\otimes p}$, and with the evident
modification of the functor $(\slot)^{\phi}$.
\end{rem}

\begin{eg}
\label{eg:18}
The coalgebra $\A^{1}(k)\in \coalg\grvect$ 
has basis 
\[
\{\beta_{0}(k), \beta_{1}(k), \dots \mid |\beta_{i}(k)|= i \,k \} 
\]
and coproduct given by 
\[
\Delta(\beta_{n}(k)) = \sum_{i+j=n} \beta_{i}(k)\otimes \beta_{j}(k).
\]
One has $\A^{1}(k)^{\phi} = \A^{1}(2k)$ and the Verschiebung map is
given by 
\begin{align*}
\beta_{2n}(k) &\mapsto \beta_{n}(2k) \\
\beta_{2n+1}(2k) &\mapsto 0.
\end{align*}
\end{eg}

If $\hringh^{\bullet}$ is an evenly graded Hopf ring over an evenly
graded ring $\rbullet{A}$ then $(\hringh^{\bullet})^{\phi}$ is an
evenly graded Hopf ring over $\rbullet{A}$ and the Verschiebung map
\[
\hringh^{\bullet} \to (\hringh^{\bullet})^{\phi}
\]
is a map of Hopf rings over $\rbullet{A}$.  For a coalgebra $C$, the
Verschiebung provides two maps 
map of evenly graded $\rbullet{A}$-algebras
\begin{align*}
\hringh^{\bullet}(C^{\phi}) &\to \hringh^{\bullet}(C) \\
\hringh^{\bullet}(C) &\to(\hringh^\bullet)^{\phi}(C)
\end{align*}
related by the commutative diagram 
\[
\xymatrix{
\hringh^{\bullet}(C^{\phi})  \ar[r]\ar[d]  &
\hringh^{\bullet}(C) 
\ar[d] \\
(\hringh^\bullet)^{\phi}(C^{\phi}) \ar[r]
&(\hringh^\bullet)^{\phi}(C)\mathrlap{\ .}
}
\]
\begin{rem}
The grading that is doubled by $\phi$ is the internal grading in
$\grvect$ and is independent of the ``even grading'' $(\slot)^{\bullet}$.
\end{rem}

\begin{prop}
\label{thm:10}
The Verschiebung map 
\begin{equation}
\label{eq:10}
\Phi^{k}(\hringmurw^{\bullet}) \to \Phi^{2k}\hringmurw^{\bullet}
\end{equation}
is a surjective map of evenly graded objects of $\grvect$.
\end{prop}

\begin{pf}
Here is the idea of the proof.  In the framework used in~\cite{RW:HR}
the Hopf ring $\Phi^{k}(\hringmurw^{\bullet})$ is ``generated'' by the
$b_{i}(k)$ and the constants $[a]\in \rbullet{MU}$.  Since, by
Example~\ref{eg:18}, the generators are all in the image
of~\eqref{eq:10}, the map must be surjective.

To make this precise we will make use of material
in~\cite{hill18:_univer_hopf}.   Since the structure of an evenly
graded Hopf ring over $\rbullet{MU}$ is given by a collection of maps between
products, for any symmetric monoidal category $\cat C$, the forgetful functor
\begin{equation}
\label{eq:11}
\hopfRingsMU(\cat C) \to(\coalg\cat C)^{\bullet}
\end{equation}
commutes with reflexive coequalizers and therefore sends effective
epimorphism to effective epimorphisms.  So while the map of coalgebras
underlying a surjective map in $\hopfRingsMU(\grvect)$ need not be
surjective in general, it is if it is part of a reflexive coequalizer
diagram.  If $\cat C$ has enough colimits then the forgetful
functor~\eqref{eq:11} has a left adjoint which we denote
$\freeHopfRing$.  The Hopf ring $\Phi^{k}\hringmurw^{\bullet}$ is
constructed from the (reflexive) coequalizer diagram
\[
\freeHopfRing(\A^{2}(k)\amalg \A^{1}(k)\amalg \A^{0})
\rightrightarrows \freeHopfRing(\A^{1}(k)) \to
\Phi^{k}\hringmurw^{\bullet}
\]
in which $\A^{i}(k)$ is regarded as an evenly graded object
concentrated entirely in degree $2$, and the parallel arrows express
the identities of Definition~\ref{def:5} on the summand
corresponding to $\A^{2}(k)\amalg \A^{0}$ and are the identity on $\A^{1}(k)$.
This implies that the map 
\[
\freeHopfRing(\A^{1}(k)) \to \Phi^{k}\hringmurw^{\bullet}
\]
is an effective epimorphism, so to show that the Verschiebung 
\[
\Phi^{k}\hringmurw^{\bullet}\to \Phi^{2k}\hringmurw^{\bullet}
\]
is an effective epimorphism it suffices to show that 
\[
\freeHopfRing(\versch):\freeHopfRing(\A^{1}(k))\to 
\freeHopfRing(\A^{1}(2k))
\]
is one.  Since $\freeHopfRing$ is a left adjoint this reduces to
showing that the Verschiebung map 
\[
\A^{1}(k)\to \A^{1}(2k)
\]
is an effective epimorphism.   This is easily checked directly and
a formula can be derived by taking the linear dual of the reflexive
equalizer of graded $\ftwo$ algebras 
\[
\ftwo[t] \rightrightarrows \ftwo[x,y] \to \ftwo[x,y]/(x^{2}-y^{2}),
\]
in which the horizontal arrows send $t$ to $x^{2}$ and $y^{2}$ respectively.
\end{pf}

\subsection{The Verschiebung ideal}
\label{sec:our-main-decomp}

Suppose that $C\in\coalg\grvect$ is a coalgebra equipped with a point $x:\unit\to
C$.   From this data one can construct a coalgebra $K$ by the pullback
square 
\[
\xymatrix{
K  \ar[r]\ar[d]  & C
\ar[d]^-{\versch} \\
\unit  \ar[r]_-{x}        & C^{\phi}\mathrlap{\ .}
}
\]
When $C$ is a Hopf algebra and $x=[0]$ is the additive unit, then the
diagram above is a pullback diagram of abelian groups in
$\coalg\grvect$ and $K$ is the {\em Verschiebung kernel}.  When $C$ is
a Hopf ring then $K$ is an ideal we will call the {\em Verschiebung
ideal}.  If $C=\hringh^{\bullet}$ is an evenly graded Hopf ring over
an evenly graded ring $A^{\bullet}$ then the Verschiebung ideal
$K=K^{\bullet}$ is an evenly graded module over $\hringh^{\bullet}$
and the map $K^{\bullet}\to \hringh^{\bullet}$ is a module map.  We
are interested in one particular example of this construction.

\begin{defin}
\label{def:13} The {\em (Ravenel-Wilson) Verschiebung ideal}
$\rwHringIdeal^{\bullet}$ is defined by the pullback square
\[
\xymatrix{
\rwHringIdeal^{\bullet}  \ar[r]\ar[d]  & \Phi^{1}\hringmurw^{\bullet}
\ar[d] \\
\unit  \ar[r]_-{[0]}        & \Phi^{2}\hringmurw^{\bullet}\mathrlap{\ .}
}
\]
\end{defin}

We now turn to our main decomposition.

\begin{prop}
\label{thm:11}
For each $k$ the Hopf algebra $\hringmurw^{2k}$ is spacelike.
\end{prop}

\begin{pf}
This is part of Proposition~\ref{thm:3}.
\end{pf}

Combined with Proposition~\ref{thm:11}, Proposition~\ref{thm:7} gives the decomposition
\begin{equation}
\label{eq:12}
\hringmurw^{2k} \approx (\hringmurw^{2k})'\otimes \Z[MU^{2k}].
\end{equation}

Since the Verschiebung is a bijection on group like elements, the
Ravenel-Wilson Verschiebung ideal could just as well be defined by the
pullback square
\begin{equation}
\label{eq:13}
\xymatrix{
\rwHringIdeal^{2k}  \ar[r]\ar[d]  & \Phi^{1}(\hringmurw^{2k})'
\ar[d] \\
\unit  \ar[r]_-{[0]}        & \Phi^{2}(\hringmurw^{2k})'\mathrlap{\ .}
}
\end{equation}
Now the two right terms in~\eqref{eq:13} are free commutative algebras
by Theorem~\ref{thm:3} and the right vertical map is surjective by
Proposition~\ref{thm:10}.   This implies

\begin{prop}
\label{thm:12}
For each $k$ there there exists an commutative
algebra section of the Verschiebung map in~\eqref{eq:13}.   Associated
to a choice of section is an isomorphism of commutative rings 
\begin{equation}
\label{eq:14}
\Phi^{1}\hringmurw^{2k} \approx K^{2k}\otimes
(\Phi^{2}\hringmurw^{2k})'\otimes \ftwo[MU^{2k}].
\end{equation}\qed
\end{prop}

\subsection{A further technical result}
\label{sec:furth-tech-result}

Our main computation requires Proposition~\ref{thm:13} below, which is
an additional technical result of Ravenel-Wilson~\cite{RW:HR}.  It is
not stated explicitly in~\cite{RW:HR} as a theorem, but appears in an
inductive argument as~\cite[(4.19)]{RW:HR} where it is a composition
of isomorphisms~\cite[(4.19)]{RW:HR} and~\cite[(4.20)]{RW:HR}, and the
statement immediately following~\cite[(4.20)]{RW:HR}.

Suppose that $k$ is a commutative ring and that $\hringh^{\bullet}$ an
evenly graded Hopf ring in $\gmodules{k}$.  The graded abelian group
structure makes each $\hringh^{2n}$ into a Hopf algebra.  As in 
Example~\ref{eg:3} we denote the coproduct and unit by
\begin{align*}
\psi &:\hringh^{2n} \to \hringh^{2n}\otimes \hringh^{2n} \\
\epsilon &:\hringh^{2n} \to k,
\end{align*}
There is a map
\[
\circ :\hringh^{2n}\otimes \hringh^{2m}\to \hringh^{2n+2m}
\]
corresponding to multiplication, and there are constants
\begin{align*}
[0]=[0]_{2n} & \in \hringh^{2n} \\
[1] &\in \hringh^{0}.
\end{align*}
The multiplicative unit for the Hopf algebra structure is the constant
$[0]$, so we have
\[
1=[0]\in \hringh^{2n}.
\]   

The kernel of $\epsilon$ is the
{\em augmentation ideal} $I^{2n}\subset \hringh^{2n}$ and the $k$-module
\begin{equation}
\label{eq:63}
Q\hringh^{2n} = I^{2n}/(I^{2n})^{2}
\end{equation}
is the module of indecomposables.   The $\circ$ product makes
$Q\hringh^{\bullet}$ into an evenly graded commutative ring.  

The fact that multiplication by $0$ in a ring is $0$ is expressed in a
Hopf ring by the identity 
\[
1\circ r = [0]\circ r= \epsilon(r).
\]
If $e\in \hringh^{2m}$ is primitive and $a,b\in I^{2n}$ then
\[
e\circ(a b) = \epsilon(a)(e\circ b)+(e\circ a)\epsilon(b) = 0
\]
so that ``circle product with $e$'' gives a map 
\[
e\circ(\slot):Q\hringh^{2n} \to \hringh^{2m+2n},
\]
whose image lies in the group of primitive elements of $\hringh^{2m+2n}$.

Let 
\[
b_{i}=x_{\ast}(\beta_{i})
\]
be the image of $\beta_{i}$ under the additive curve 
\[
x:\addGrp^{1}\to H_{\ast}(\umu^{2};\ftwo).
\]
From the assumption
$\beta_{0}^{\ast}(x)=[0]$ we have $b_{0}=1$, and the element $b_{1}$
is primitive.

\begin{prop}[\cite{RW:HR}]
\label{thm:13}
The map 
\[
b_{1}\circ(\slot):QH_{\ast}(\umu^{2k};\ftwo) \to H_{\ast}((\umu^{2k+2})';\ftwo)
\]
is an isomorphism of $QH_{\ast}(\umu^{2k};\ftwo)$ with the primitives in
\[
H_{\ast}((\umu^{2k+2})';\ftwo),
\]
where $(\umu^{2k+2})'\subset \umu^{2k+2}$ is the connected component
containing the base point.   In particular it is a monomorphism.
\end{prop}

By Theorem~\ref{thm:5}, the above has a purely algebraic counterpart.

\begin{cor}
\label{thm:14}
The map 
\[
b_{1}\circ:Q(\hringmur^{2k}\otimes \ft)\to \hringmurw^{2k+2}\otimes \ft
\]
induces an isomorphism of $Q(\hringmur\otimes\ft)^{2k}$ with the primitives in
\[
(\hringmurw^{2k+2}\otimes\ft)'.
\]
In particular it is a monomorphism.
\end{cor}

\begin{rem}
\label{rem:6}
Since primitive elements are annihilated by the Verschiebung, the
image of $b_{1}\circ(\slot)$ is contained in the Verschiebung ideal.
This also follows from the fact that $\versch(b_{1}(1))=0$
(Example~\ref{eg:18}) since the Verschiebung ideal is an ideal.
\end{rem}

For more on the definition and summary of the useful formulae for
graded Hopf rings in the category of graded abelian groups the reader
is referred to~\cite[Lemma~1.12]{RW:HR}.

\section{The equivariant theory}
\label{sec:equivariant-theory}

In this section we turn to equivariant stable homotopy theory and the
Hopf rings arising from real oriented cohomology theories.  Our main
purpose is to set up the equivalence $\mot$ between the category of
evenly graded free abelian groups and the category of {\em pure
modules} over the equivariant Eilenberg-MacLane spectrum $H\zm$.
Using this we transport the Ravenel-Wilson Hopf ring to equivariant
homotopy theory and, appealing to Theorem~\ref{thm:6}, construct the
map~\eqref{eq:20} appearing in the statement of Theorem~\ref{thm:26},
our main result.

General references for equivariant stable homotopy theory
\cite{MR1413302,MR1361893,LMayS,MR1922205,schwede:_lectur},
and~\cite[\S2 and Appendix~B]{MR3505179}

\subsection{Equivariant $R$-modules}
\label{sec:equiv-r-modul}

For a $G$-equivariant $E_{\infty}$ ring $R$ let $\rmod{R}$ denote the
homotopy category of equivariant left $R$-modules and equivariant
$R$-module maps, regarded as an additive category. The category
$\rmod{R}$ has arbitrary coproducts and products.  It becomes an
strongly additive symmetric monoidal category under 
\[
M\otimes N = M\underset{R}{\wedge}N,
\]
with $R$ as the tensor unit.  The functor 
\[
\freeRmodule{R}{X} = R\wedge X_{+}
\]
from the homotopy category $\ho\spaces^{G}$ of $G$-spaces and
equivariant maps to $\rmod{R}$ is symmetric monoidal when the category
of pointed $G$-spaces is equipped with the symmetric monoidal
structure given by the Cartesian product.

\subsection{Restriction and geometric fixed points}
\label{sec:restr-geom-fixed}

Associated to an inclusion $i:H\subset G$ there is a strongly additive symmetric
monoidal restriction functor
\[
\res_{i}= \res:\rmod{R} \to \rmod{\res R}.
\]
The functor
\begin{align*}
\rmod{\res R} & \to \rmod{R} \\
X &\mapsto G_{+}\underset{H}{\wedge} X
\end{align*}
is both a left and right adjoint to $\res$ (the Wirthm\"uller
isomorphism).  There is also a symmetric monoidal {\em geometric fixed
point} functor~\cite[\S B.10.6]{MR3505179}
\[
\gfp{G}:\rmod{R}\to \rmod{\gfp{G}R}.
\]

\begin{prop}
\label{thm:15}
A map $X\to Y$ in $\rmod{R}$ is an isomorphism if and only if for
every $H\subset G$ the map 
\[
\gfp{H}X\to \gfp{H}Y
\]
is an isomorphism.
\end{prop}

\begin{pf}
This follows directly from, say, \cite[Proposition~2.52]{MR3505179}.
\end{pf}

The following result is useful for analyzing the symmetric monoidal
structure on $\rmod{R}$.

\begin{lem}
\label{thm:16} Suppose that $G$ is a group and $E$ is a $G$-equivariant
homotopy commutative ring with the property that the map
\[
\pi_{0}^{G}(S^{0}) \to \pi_{0}^{G}E
\]
induced by the unit $S^{0}\to E$ factors through the restriction
\[
\pi_{0}^{G}S^{0}\to \pi_{0}S^{0} =  \Z.
\]
If $V$ is a representation of $G$ then the symmetry map $S^{V}\wedge
S^{V}\to S^{V}\wedge S^{V}$ induces
\[
(-1)^{\dim V}: E\wedge S^{V}\wedge S^{V}\to E\wedge S^{V}\wedge S^{V}.
\]
\end{lem}

\begin{pf}
The symmetry map
\[
S^{V}\wedge S^{V}\to S^{V}\wedge S^{V}
\]
defines an element 
\[
\epsilon\in \pi_{2V}^{G}(S^{2V}) \approx \pi_{0}^{G}(S^{0}),
\]
and the assertion is that, under the map induced by the unit, the
image of this element in $\pi_{0}^{G}E$ is $(-1)^{\dim V}$.  By
assumption, the map
\[
\pi_{0}^{G}S^{0}\to \pi_{0}^{G}E
\]
factors through the restriction 
\[
\pi_{0}^{G}S^{0}\to \pi_{0}S^{0}=\Z.
\]
The result follows from the fact that the degree
of the underlying non-equivariant map of spheres is $(-1)^{\dim V}$.
\end{pf}

\begin{rem}
\label{rem:7} Lemma~\ref{thm:16} does not hold without smashing with
$E$.  The symmetry map corresponds to the element of the
Burnside ring mapping to $(-1)^{\dim V^{H}}$ under the character
sending a $G$-set to its $H$ fixed points.  For instance when $G=\zt$
and $V=\rho$ this element is $(-1+[\zt])$.
\end{rem}

\subsection{Pure modules}
\label{sec:pure-modules}

We now specialize to the case $G=\zt$.  Following~\cite{MR3505179} we
let $\sigma$ be the sign representation and $\rho$ the real regular
representation of $G$.  Denote by $\Zm$ the constant Mackey
functor $\Z$, and $H\Zm$ the corresponding Eilenberg-MacLane spectrum.

\subsubsection{The category of pure modules}
\label{sec:categ-pure-modul}

Let $\hzmodf\subset \hzmod$ be the smallest full subcategory of
$\hzmod$ containing the objects $H\zm\wedge S^{2k\rho}$ for all
$\kk$, and closed under arbitrary coproducts.  The subcategory $\hzmodf$ is
a strongly additive symmetric monoidal category.  We will call objects of
$\hzmodf$ {\em pure}.

\begin{lem}
\label{thm:17}
For integers $k, \ell \in \Z$, the restriction mappings
\begin{align*}
\hzmodmap{H\Zm\wedge S^{k\rho}}{H\Zm\wedge S^{\ell\rho}} &\to
\zmodmap{H\Z\wedge S^{2k}}{H\Z\wedge S^{2\ell}} \\
\hzmodmap{H\zm\wedge S^{k\rho-1}}{H\Zm\wedge S^{\ell\rho}} &\to
\zmodmap{H\Z\wedge S^{2k-1}}{H\Z\wedge S^{2\ell}} \\
\end{align*}
are isomorphisms, and so
\[
\hzmodmap{H\Zm\wedge
S^{k\rho}}{H\Zm\wedge S^{\ell\rho}}  =
\begin{cases}
\Z & k=\ell\\
0 & k\ne \ell,
\end{cases}
\]
and for all $k$ and $\ell$
\[
\hzmodmap{H\Zm\wedge
S^{k\rho-1}}{H\Zm\wedge S^{\ell\rho}} =0.
\]
\end{lem}

\begin{pf}
In the language of~\cite{MR3505179} this is just the assertion that
$H\Zm\wedge S^{k\rho}$ is a $2k$-slice and $H\Zm\wedge S^{k\rho-1}$
is a $(2k-1)$-slice.   It is easily checked directly.   By tensoring
with the inverse of $H\zm\wedge S^{k\rho}$ and using the isomorphism 
\[
\hzmodmap{H\zm}{H\zm\wedge X} = H_{0}(X;\zm) = H^{0}(DX;\zm)
\]
one reduces to checking
that the restriction mapping 
\[
\tilde H^{0}_{\zt}(S^{0},\zm)\to \tilde H^{0}(S^{0};\Z)
\]
is an isomorphism and that for all $\ell>0$
\[
\tilde H^{\zt}_{0}(S^{\ell \rho});\Zm) = \tilde H^{0}_{\zt}(S^{\ell\rho};\zm)
= 0.
\]
The first is the definition of the constant Mackey functor and the
second follows from the fact that $S^{\ell\rho}$ is obtained from
$S^{\ell}$ by attaching free cells of dimension $(\ell+1)$ and
higher (see \cite[\S3.3]{MR3505179}).
\end{pf}

\begin{cor}
\label{thm:18}
If $X$ and $Y$ are pure $H\zm$-modules then the restriction mapping 
\[
\hzmodmap{X}{Y} \to 
\zmodmap{\res X}{\res Y}
\]
is an isomorphism. \qed
\end{cor}

\subsubsection{Pure modules and graded free abelian groups}
\label{sec:pure-modules-graded}

There is a unique coproduct preserving functor
\begin{equation}
\label{eq:15}
\mot:\frab\to \hzmodf
\end{equation}
from the category of evenly graded free abelian groups, satisfying
\[
\mot(\Z[2n])= S^{n\rho}\wedge H\zm.
\]
Propositions~\ref{thm:16} and~\ref{thm:17} imply

\begin{prop}
\label{thm:19}
The functor~\eqref{eq:15} is an equivalence of strongly additive symmetric monoidal
categories, with inverse given by $\pi_{\bullet}\res$.   
\qed
\end{prop}

\subsubsection{Projective space}
\label{sec:projective-space}

For $n\le \infty$ we make $\cp^{n}$ into a $\zt$-space by identifying it with the
space of complex lines in $\C^{n+1}$ and letting $\zt$ act by complex
conjugation.   The quotient space
\[
\cp^{n}/\cp^{n-1}
\]
is identified with the one point compactification of $\C^{n}$ which,
in turn is canonically isomorphic to $S^{n\rho}$.   There is thus a
cofibration sequence
\begin{equation}
\label{eq:16}
\cp^{n-1}\to \cp^{n}\to S^{n\rho}.
\end{equation}
Using this and Lemma~\ref{thm:17}, one can easily establish the
following result

\begin{prop}
\label{thm:20}
For all $0\le n_{1},\dots,
n_{k}\le \infty$ the $H\zm$ module 
\[
\freeRmodule{H\Zm}{\cp^{n_{1}}\times
\cdots \times \cp^{n_{k}}}
\]
is pure. \qed
\end{prop}

\begin{cor}
\label{thm:21} For each $k\ge 0$ there is a  unique coalgebra
isomorphism
\begin{equation}
\label{eq:17}
\mot(\A^{k}) \approx \freeRmodule{H\zm}{(\cp^{\infty})^{k}}
\end{equation}
restricting to the isomorphism
\[
\motZ(\A^{k}) = \freeRmodule{H\Z}{(\cp^{\infty})^{k}} = H\Z\wedge (\cp^{\infty})^{k}_{+}
\]
of Remark~\ref{rem:3}.  Under this equivalence, the
isomorphism~\eqref{eq:17} is an isomorphism of Hopf algebras. \qed
\end{cor}

\subsection{Grading}
\label{sec:grading}

In equivariant homotopy theory it is useful to grade things over the
real representation ring $RO(G)$.  We will follow the convention,
introduced by Hu and Kriz in~\cite{MR1808224}, of using the symbol
$\star$ as a wildcard matching a real virtual representation, and
using $M_{\star}$ and $M^{\star}$ to denote $RO(G)$-graded objects.

Associated to  an $RO(G)$-graded object $M^{\star}$ and virtual
representation $V\in RO(G)$ of virtual dimension $d$ is the 
the $d\Z$-graded object $M(V)^{dk} = M^{kV}$.   This construction 
is compatible with the restriction map along the inclusion of a
subgroup $H\subset G$.

We are interested in the special case $G=\zt$ and with $V=\rho$
the real regular representation.  In this case $M(\rho)^{\bullet}$ is
an evenly graded object, with $M(\rho)^{2n}=M^{n\rho}$. 

\begin{eg}
\label{eg:19}
If $E$ is a $\zt$-spectrum, then $\underline{E}(\rho)^{\bullet}$ is
the evenly graded $\zt$-space with 
\[
\underline{E}(\rho)^{2n} = \underline{E}^{n\rho} =
\Omega^{\infty}S^{n\rho}\wedge E.
\]
In this case 
\[
\res \underline{E}(\rho)^{2n} = \underline{E}^{2n}.
\]
\end{eg}

\begin{eg}
\label{eg:20}
If $E$ is a $\zt$-spectrum then $E(\rho)^{\bullet}$ is the evenly
graded abelian group with 
\[
E(\rho)^{2n}= E^{n\rho}_{\zt}(\text{pt}) =\pi_{-n\rho}^{\zt}E.
\]
\end{eg}

\subsection{Real oriented spectra and real bordism}

\subsubsection{Real orientations}
\label{sec:real-orientations}

\begin{defin}
\label{def:14}
A {\em real oriented} spectrum is a $\zt$-equivariant homotopy
commutative ring spectrum $E$, equipped with an element $x\in
\tilde E_{\zt}^{\rho}(\cp^{\infty})$ whose restriction to $\tilde E_{\zt}^{\rho}(\cp^{1})$
corresponds to the element $1$ under the isomorphism 
\[
\tilde E_{\zt}^{\rho}(\cp^{1}) \approx \tilde E_{\zt}^{\rho}(S^{\rho}) \approx
E^{0}(\text{pt}). 
\]
\end{defin}

From the cofibration sequence~\eqref{eq:16} one can easily prove

\begin{prop}[Araki\cite{MR614829}, Landweber\cite{MR0222890}]
If $E$ is a real oriented spectrum then the maps
\begin{align*}
E^{\zt}_{\star}\LL x\RR &\to E_{\zt}^{\star}(\cp^{\infty}) \\
E^{\zt}_{\star}\LL x,y\RR &\to E_{\zt}^{\star}(\cp^{\infty}\times \cp^{\infty})
\end{align*}
are isomorphisms, where, in the bottom expression, $x$ and $y$ are the
pullbacks of the real orientation $x$ along the projections to the
first and second factors.  
\end{prop}
The equivariant $E$-cohomology homomorphism
induced by map
\[
\cp^{\infty}\times\cp^{\infty}\to \cp^{\infty}
\]
classifying the tensor product of the two universal ``real'' line
bundles sends the real orientation $x$ to 
\[
F(x,y) = \sum a_{ij}\, x^{i}y^{j} \in 
E_{\zt}^{\star}(\cp^{\infty}\times \cp^{\infty}) \approx E^{\zt}_{\star}\LL x,y\RR 
\]
with
\[
a_{ij}\in E^{\zt}_{(i+j-1)\rho} = E(\rho)^{-(i+j-1)\rho}.
\]
As in the non-equivariant case, this defines a formal group law over
$E_{\zt}^{\star}$, and so, with the conventions of
\S\ref{sec:grading}, a map of evenly graded rings
\[
\rbullet{MU} \to E(\rho)^{\bullet}.
\]

\subsubsection{The real bordism spectrum}

The universal example of a real oriented spectrum is the {\em real
bordism} spectrum $\mur$ of Landweber~\cite{MR0222890} and
Fujii~\cite{MR0420597}, and later investigated by
Araki~\cite{MR614829}, Hu-Kriz~\cite{MR1808224}, and
in~\cite{MR3505179}.    The spectrum $\mur$ is the $\zt$ equivariant
spectrum, constructed as Thom spectrum of the universal complex vector bundle over
$BU$, with $\zt$ acting by complex conjugation.   It can be realized
as the $\zt$-equivariant spectrum 
\begin{equation}
\label{eq:58}
\varinjlim S^{-n\rho}\wedge MU_{\R}(n)
\end{equation}
associated to the maps 
\begin{equation}
\label{eq:57}
S^{\rho}\wedge MU_{\R}(n-1)\to MU_{\R}(n)
\end{equation}
in which $MU_{\R}(m)$ denotes the Thom complex of the universal complex
$m$-plane bundle over $BU(m)$ with $\zt$ acting by complex
conjugation.   The inclusion of the zero section
\[
\cp^{\infty}\to MU_{\R}(1)
\]
is an equivariant weak equivalence, and from~\eqref{eq:58} gives the
universal real orientation 
\[
\Sigma^{\infty}\cp^{\infty}\to
S^{\rho}\wedge \mur^{\rho}
\]
of $\mur$.

The Schubert cell decomposition of Grassmannians equips the spectrum
$\mur$ with a stable cell decomposition into cells of the form
$D(n\rho)$.   As observed by Araki~\cite{MR614829} this implies

\begin{prop}
\label{thm:22}
For every $n\in\Z$ the spectrum 
\[
H\zm\wedge S^{n\rho}\wedge\mur
\]
is pure.   \qed
\end{prop}

A much deeper result holds.

\begin{prop}[Araki~{\cite[Theorem~4.6]{araki78:_coeff_mr}},
Hu-Kriz~{\cite[Theorem~2.28]{MR1808224}}, {\cite[Theorem~6.1]{MR3505179}}]
\label{thm:24}
The map
\[
\rbullet{MU} \to \mur(\rho)^{\bullet}
\]
is an isomorphism, with inverse given by the restriction mapping
\[
\res:\pi_{k\rho}^{\zt}\mur \to \pi_{2k}MU.
\]
\end{prop}

As in~\cite[Lemma~2.17]{MR1808224}, one simple consequence of Proposition~\ref{thm:24} is

\begin{cor}
\label{thm:25}
If $E$ is real oriented then $E$ satisfies the condition of
Lemma~\ref{thm:16} and so the ring $E(\rho)^{\bullet}$ is an evenly
graded commutative ring. \qed
\end{cor}

\subsection{Hopf rings from real oriented spectra}
\label{sec:real-orient-spectra}

Suppose that $E$ is a real oriented spectrum.  As in the previous
section let $\rbullet{E(\rho)}$ be the evenly graded ring given by
\[
E(\rho)^{2n} = E^{n\rho}_{\zt}(\text{pt}),
\]
and
\begin{equation}
\label{eq:18}
\rbullet{MU}\to \rbullet{E(\rho)}
\end{equation}
the map classifying the formal group law.

By Corollary~\ref{thm:25} and Example~\ref{eg:19} the spaces 
\[
\underline{E}(\rho)^{\bullet}
\]
form an evenly graded Hopf ring in $\spaces^{\zt}$ over the evenly
graded commutative ring $E(\rho)^{\bullet}$.   Restricting scalars
along~\eqref{eq:18} makes
$\underline{E}(\rho)^{\bullet}$ into an evenly graded Hopf ring in $\spaces^{\zt}$
over $\rbullet{MU}$.    
As in \S\ref{sec:exampl-from-compl} this equips
\[
\freeRmodule{H\Zm}{\underline{E}(\rho)^{\bullet}}= H\Zm\wedge \underline{E}(\rho)^{\bullet}_{+}
\]
with an additive curve
\[
x:\mot(\A^{1})=\freeRmodule{H\Zm}{\cp^{\infty}} \to \freeRmodule{H\Zm}{\underline{E}(\rho)^{2}}
\]
which, by Theorem~\ref{thm:6}, is classified by a unique map
\begin{equation}
\label{eq:19}
\mot(\hringmurw^{\bullet}) \to  \freeRmodule{H\Zm}{\underline{E}(\rho)^{\bullet}}
\end{equation}
of evenly graded Hopf rings over $\rbullet{MU}$.

This applies in particular to the case $E=\mur$.   For simplicity
write 
\[
\hringmur^{\bullet} = \freeRmodule{H\zm}{\underline{\mur}(\rho)^{\bullet}}.
\]
Theorem~\ref{thm:2}, our main result, is a consequence of the more refined

\begin{thm}
\label{thm:26}
The induced map
\begin{equation}
\label{eq:20}
\mot(\hringmurw^{\bullet}) \to \hringmur^{\bullet}
\end{equation}
is an isomorphism.
\end{thm}

Note that Theorem~\ref{thm:5} implies

\begin{prop}
\label{thm:27}
The restriction
\[
\res\mot(\hringmurw^{\bullet}) \to \res \hringmur^{\evenstar}
\]
of~\eqref{eq:20} is an isomorphism. \qed
\end{prop}

By Proposition~\ref{thm:15} this means that to prove
Theorem~\ref{thm:26} it suffices to show that~\eqref{eq:20} induces an
isomorphism of geometric fixed points.  

\subsection{Modified geometric fixed points}
\label{sec:geom-fixed-points}

The discussion  will take place in more ordinary terms if we slightly modify
the geometric fixed point functor
\begin{equation}
\label{eq:21}
\gfp{\zt}:\hzmod \to \modules_{\gfp{\zt} H\zm}.
\end{equation}
Since
\[
\phig(H\zm) = \hft[u]\qquad |u|=2,
\]
(see for example~\cite[Proposition~7.5]{MR3505179}), the
functor~\eqref{eq:21} is a strongly additive symmetric monoidal functor
\[
\phig:\hzmod \to \hztumod.
\]
This can be composed with 
\[
X \mapsto \hft\underset{\hft[u]}{\wedge}X
\]
to give a further strongly additive symmetric monoidal functor 
\[
\bphig:\hzmod \to \hztmod.
\]

\begin{eg}
\label{eg:21}
If $Y$ is a $\zt$ space then 
\[
\bphig(H\Zm\wedge Y_{+}) = H\ftwo\wedge Y^{\zt}_{+}.
\]
\end{eg}

\begin{defin}
\label{def:15}
A $\zt$ spectrum $X$ is {\em bounded below} if the fixed point
spectrum $X^{\zt}$ and the underlying spectrum $\res X$ have the
properties
\begin{align*}
\pi_{i} X^{\zt} &=0 \qquad i \ll 0 \\
\pi_{i} \res X &=0 \qquad i \ll 0 \\
\end{align*}
\end{defin}

\begin{lem}
\label{thm:28}
A map $X\to Y$ in $\hztumod$ of objects which are bounded below is an
isomorphism if and only if 
\[
\hft\underset{\hft[u]}{\wedge}X\to
\hft\underset{\hft[u]}{\wedge}Y
\]
is an isomorphism in $\hztmod$.
\end{lem}

\begin{pf}
By passing to the mapping cone we reduce to the statement that if $X$
is bounded below and $\hft\underset{\hft[u]}{\wedge}X$ is contractible
then $X$ is contractible.   This follows easily from the cofibration sequence
\[
S^{2}\wedge X \xrightarrow{u}{} X\to 
\hft\underset{\hft[u]}{\wedge}X.
\]
\end{pf}

Since geometric fixed points preserves connectivity,
Lemma~\ref{thm:28} and Proposition~\ref{thm:15} give

\begin{cor}
\label{thm:29}
A map $X\to Y$ in $\hzmod$ of $H\Zm$ modules which are bounded below
is an isomorphism if and only if the maps
\begin{align*}
\res X &\to \res Y \\
\bphig(X) &\to \bphig(Y)
\end{align*}
of underlying and modified geometric fixed points are isomorphisms. \qed
\end{cor}

Since the spectra $\mot(\hringmurw^{2m})$ and $\hringmur_{2m}$ are
bounded below, Proposition~\ref{thm:27} and Corollary~\ref{thm:29}
reduce Theorem~\ref{thm:26} to showing that the map
\begin{equation}
\label{eq:22}
\bphig\mot(\hringmurw^{\bullet}) \to \bphig(\hringmur^{\bullet})
\end{equation}
is an isomorphism.   
\section{The Hopf ring $X^{\bullet}$}
\label{sec:hopf-ring-xbullet}
The map~\eqref{eq:22} is a map in the category $\hztmod$ which, by
Proposition~\ref{thm:8} is equivalent to $\grvect$ via the functors
\[
\vmot:\grvect \leftrightarrow \hztmod:\pi_{\ast}.
\]  
To locate
everything in $\grvect$ first note that the following diagram commutes
\begin{equation}
\label{eq:23}
\xymatrix{
\frab  \ar[r]^-{\mot}\ar[d]_{\Phi^{1}}  & \hzmod
\ar[d]^{\bphig } \\
\grvect  \ar[r]_-{\vmot}        & \hztmod\mathrlap{,}
}
\end{equation}
where $\Phi^{1}$ is the functor defined in~\S\ref{sec:some-algebra-i}.
This means that we may identify the image of left hand side
of~\eqref{eq:22} in $\grvect$ with
\[
\Phi^{1}(\hringmurw^{\bullet}).
\]
The image of the right hand side is a Hopf ring over $MU^{\bullet}$ 
equipped with the additive curve 
\[
u=\pi_{\ast}\bphig(x):\A^{1}(1)=\pi_{\ast}\bphig \mot(\A^{1}) \to
\pi_{\ast}\bphig(\hringmur^{\bullet}), 
\]
and the image of the map~\eqref{eq:22} is the unique map of evenly
graded Hopf rings over $MU^{\bullet}$ classifying this additive curve.

The purpose of this section is to identify the right hand side
of~\eqref{eq:22} with the mod $2$ homology of an evenly graded Hopf
ring $X^{\bullet}$ over $MU^{\bullet}$ in spaces, and from that,
establish the decompositions~\eqref{eq:59} and~\eqref{eq:60}, their
compatiblity, and the commutativity of the diagrams~\eqref{eq:34} and
~\eqref{eq:35}.

In all of the remaining discussion we will be concerned
solely with homology with coefficients in $\ft$.  To simplify the
notation we will now use the symbol $H_{\ast}(\slot)$ to mean
$H_{\ast}(\slot;\ft)$.   

\subsection{The spaces $X^{2k}$}
\label{sec:spaces-x2k}

The relationship between geometric fixed
points and suspension spectra (Example~\ref{eg:21}) gives 
\begin{align*}
\bphig\hringmur^{2k} &=
\bphig\freeRmodule{H\zm}{\umur(\rho)^{2k}} \\
&=
\bphig (H\zm\wedge \umur^{k\rho}{}_{+})\\
&= H\ftwo\wedge (\umur^{k\rho})^{\zt}_{+}.
\end{align*}
Let 
\[
X^{2k} = (\umur^{k\rho})^{\zt}.
\]
By functoriality, the sequence of spaces $X^{2k}$ forms an evenly graded
Hopf ring $X^{\bullet}$ over $\rbullet{MU}$ in $\ho\spaces$.
The constant associated to $[a]\in MU^{2k}$ is obtained by passing to
fixed points from the unique equivariant map 
\[
S^{0}\to \umur^{k\rho}
\]
restricting to $a\in\pi_{0}\umu^{2k}=MU^{2k}(\text{pt})$.  By the
above, the Hopf ring $\pi_{\ast}\bphig \hringmur^{\bullet}$ is
$H_{\ast}(X^{\bullet};\ftwo)$, and the weight one additive curve $u$
is the homology homomorphism induced by the map
\begin{equation}
\label{eq:24}
\rp^{\infty} \to X^{2}
\end{equation}
obtained by passing to fixed points from the real orientation 
\[
\cp^{\infty}\to \umur^{\rho}.   
\]
The image of~\eqref{eq:22} in $\grvect$ is therefore the unique map 
\begin{equation}
\label{eq:25}
\Phi^{1}(\hringmurw^{\bullet}) \to H_{\ast}(X^{\bullet})
\end{equation}
of Hopf rings over $MU^{\bullet}$, classifying the weight one additive
curve coming from~\eqref{eq:24}.

\subsection{Decomposition}
\label{sec:decomposition}

Our next aim is to decompose $H_{\ast}(X^{2k})$ in a manner compatible
with the decomposition~\eqref{eq:14} of Proposition~\ref{thm:12}.
This will be achieved by exploiting the fixed point inclusion
\begin{equation}
\label{eq:26}
X^{\bullet}\to \umu^{\bullet}
\end{equation}
which is map of evenly graded Hopf rings over $MU^{\bullet}$.
Proposition~\ref{thm:24} implies that~\eqref{eq:26} induces a
bijection of path components.  By Example~\ref{eg:16} and
Proposition~\ref{thm:7} there is an algebra decomposition
\begin{equation}
\label{eq:27}
H_{\ast}(X^{2k}) \approx H_{\ast}((X^{2k})')\otimes \ftwo[MU^{2k}].  
\end{equation}

Since
\[
H_{\ast}(\umu^{\bullet}) = \Phi^{2}(\hringmurw^{\bullet})
\]
composition of~\eqref{eq:25} with the fixed point inclusion gives a
sequence
\begin{equation}
\label{eq:28}
\Phi^{1}\hringmurw^{\bullet}\to H_{\ast}(X^{\bullet})\to \Phi^{2}\hringmurw^{\bullet}
\end{equation} 
of Hopf rings over $\rbullet{{MU}}$.

\begin{prop}
\label{thm:30}
The composition~\eqref{eq:28} is the Verschiebung.
\end{prop}

\begin{pf}
This follows from Proposition~\ref{thm:9} and the fact that the fixed
point inclusion 
\[
\rp^{\infty}\to \cp^{\infty}
\]
induces the Verschiebung 
\[
\A^{1}(1) =H_{\ast}(\rp^{\infty};\ftwo) \to \A^{1}(2) = H_{\ast}(\cp^{\infty};\ftwo).
\]
\end{pf}

From Proposition~\ref{thm:10} we conclude

\begin{cor}
\label{thm:31}
The fixed point inclusion
\[
H_{\ast}(X^{\bullet};\ftwo) \to H_{\ast}(\umu^{\bullet};\ftwo) = \Phi^{2}(\hringmurw^{\bullet})
\]
is a surjective map of underlying $\ftwo$ vector spaces.  \qed
\end{cor}

To go further, smash the cofibration sequence
\[
(\zt)_{+} \to S^{0} \to S^{\sigma}, 
\]
with $S^{1}\wedge S^{(k-1)\rho}\wedge \mur$, take fixed points, and pass
to zeroth spaces.   This leads to the fibration sequence
\begin{equation}
\label{eq:29}
\xymatrix{
*+{\umu^{2k-1}}  \ar[r]  &  BX^{2(k-1)} \ar[d]   &   \\
        & X^{2k}  \ar[r]         & *+{\umu^{2k},}
}
\end{equation}
in which we have written
\[
BX^{2(k-1)}=\big(\Omega^{\infty}S^{1}\wedge S^{(k-1)\rho}\wedge \mur)^{\zt}.
\]
This is justified by
\begin{lem}
\label{thm:32}
The space $\Omega^{\infty}(S^{1}\wedge S^{(k-1)\rho}\wedge \mur)^{\zt}$ is
connected.
\end{lem}

\begin{pf}
We have 
\[
\pi_{0}\Omega^{\infty}S^{1}\wedge S^{(k-1)\rho}\wedge \mur^{\zt}
= \pi^{\zt}_{-(k-1)\rho-1}\mur
\]
which vanishes by~\cite[Theorem~1.13 and
Corollary~4.65]{MR3505179}
\end{pf}

\begin{rem}
\label{rem:8} The original proof of the vanishing result is in the
unpublished manuscript~\cite[Theorem~4.6]{araki78:_coeff_mr} of Araki.
The first published version is due to Hu and
Kriz~\cite[Theorem~4.11]{MR1808224}.
\end{rem}

By Proposition~\ref{thm:24}, the sequence~\eqref{eq:29}
restricts to a fibration sequence
\begin{equation}
\label{eq:30}
\xymatrix{
*+{\umu^{2k-1}}  \ar[r]  &  BX^{2(k-1)} \ar[d]   &   \\
        & (X^{2k})'  \ar[r]         & *+{(\umu^{2k})'\mathrlap{\ .}}
}
\end{equation}

\begin{prop}
\label{thm:33}
The sequence 
\begin{equation}
\label{eq:31}
\ft\to H_{\ast}BX^{2(k-1)} \to H_{\ast}(X^{2k})' \to
H_{\ast}(\umu^{2k})'\to \ft
\end{equation}
is a short exact sequence of Hopf-algebras.   
.
\end{prop}

\begin{pf}
By Corollary~\ref{thm:31} the map
\[
H^{\ast}((\umu^{2k})') \to H^{\ast}((X^{2k})')
\]
is a monomorphism of connected graded Hopf algebras over $\ft$ and
hence a flat map of rings (Milnor-Moore~\cite[Theorem~4.4]{MilMoore}).
It follows from the Eilenberg-Moore spectral sequence that
\[
H^{\ast}((X^{2k})')\underset{H^{\ast}((\umu^{2k})')}{\otimes}\ftwo \to H^{\ast}(BX^{2(k-1)})
\]
is an isomorphism.   The assertion to be proved is just the linear
dual of this fact.  
\end{pf}

Propositions~\ref{thm:33} and~\ref{thm:30} imply that the restriction
of the map 
\[
\Phi^{1}\hringmurw^{2k}\to H_{\ast}(X^{2k})
\]
to the Verschiebung ideal factors through $H_{\ast}(BX^{2(k-1)})$, leading
to the following diagram of short exact sequences of Hopf
algebras 
\begin{equation}
\label{eq:32}
\xymatrix{
\ft \ar[r] & \brk{2k} \ar[r]\ar[d]  &  (\Phi^{1}\hringmurw^{2k})' \ar[r]\ar[d] &H_{\ast}(\umu^{2k})'\ar[r]\ar@{=}[d]  &  \ft  \\
\ft \ar[r] & H_{\ast}BX^{2(k-1)} \ar[r]        & H_{\ast}(X^{2k})'  \ar[r]
& H_{\ast}(\umu^{2k})'\ar[r]    &\ft\mathrlap{\ .}
}
\end{equation}
As in Proposition~\ref{thm:12} one may choose a commutative algebra
section of the upper right map and, using it, construct isomorphisms
\begin{align*}
\brk{2k}\otimes H_{\ast}(\umu^{2k})' &\to (\Phi^{1}\hringmurw^{2k})' \\
H_{\ast}BX^{2(k-1)}\otimes H_{\ast}(\umu^{2k})' &\to H_{\ast}(X^{2k})' \\
\end{align*}
with respect to which the middle vertical map in~\eqref{eq:32}
is the tensor product of the left and right vertical maps.

Putting this all together we have exhibited algebra isomorphisms
\begin{gather}
\label{eq:59}
\Phi^{1}\hringmurw^{2k} \approx \brk{2k} \otimes
\Phi^{2}(\hringmurw^{2k})'\otimes \ftwo[MU^{2k}] \\
\label{eq:60}
H_{\ast}X^{2k} \approx H_{\ast}BX^{2k-2} \otimes
\Phi^{2}(\hringmurw^{2k})'\otimes \ftwo[MU^{2k}] 
\end{gather}
with respect to which the map 
\begin{equation}
\label{eq:33}
\Phi^{1}\hringmurw^{2k} \to H_{\ast}X^{2k}
\end{equation}
is the tensor product of the identity map of
$\Phi^{2}(\hringmurw^{2k})'\otimes \ftwo[MU^{2k}]$
with the map 
\[
K^{2k}\to H_{\ast}BX^{2(k-1)}.
\]

\begin{rem}
\label{rem:9}
Adjoint to the map $BX^{2(k-1)}\to X^{2k}$ is a map
$X^{2(k-1)}\to \Omega X^{2k}$.  The spectrum associated to these maps is
the geometric fixed point spectrum $\Phi^{\zt}\mur=MO$.
\end{rem}

\subsection{Changing $k$}
\label{sec:changing-k}

Recall the element
\[
\beta_{1}(1)\in \A^{1}(1)
\]
and its images (both denoted $b_{1}(1)$) in
$(\Phi^{1}\hringmurw^{2})_{1}$ and $H_{1}X^{2}$ under the universal
additive curve, and the additive curve $u$ of \S\ref{sec:spaces-x2k}.
Our proof that the map~\eqref{eq:33} is an isomorphism involves an
induction on $k$.  What relates the successive values of $k$ is the
operation $b_{1}(1)\circ(\slot)$, and the diagram
\begin{equation}
\label{eq:34}
\xymatrix{
\Phi^{1}(\hringmurw^{2(k-1)})  \ar[r]\ar[d]_{b_{1}(1)\circ(\slot)}  & H_{\ast}X^{2(k-1)}
\ar[d]^{b_{1}(1)\circ(\slot)} \\
\Phi^{1}(\hringmurw^{2k})   \ar[r]        & H_{\ast}X^{2k}\mathrlap{,}
}
\end{equation}
which commutes by virtue of the fact that
\[
\Phi^{1}\hringmurw^{\bullet}\to H_{\ast}X^{\bullet}
\]
is a map of Hopf rings.  As explained in the discussion leading up to
the statement of Proposition~\ref{thm:13} the fact that $b_{1}(1)$ is
primitive implies that the  vertical maps in~\eqref{eq:34} induce 
maps of indecomposables from the top row.   Since $b_{1}(1)$ maps to
zero in $(\Phi^{2}\hringmurw^{2})_{1}=0$, the images of the vertical
maps are contained in the Hopf algebra kernels in~\eqref{eq:32}.
This leads to the commutative diagram
\begin{equation}
\label{eq:35}
\xymatrix{
Q(\Phi^{1}(\hringmurw^{2(k-1)}))  \ar[r]\ar[d]_{b_{1}(1)\circ(\slot)}  & Q(H_{\ast}X^{2(k-1)})
\ar[d]^{b_{1}(1)\circ(\slot)} \\
\brk{2k}  \ar[r]        &
H_{\ast}BX^{2(k-1)}\mathrlap{\ .}
}
\end{equation}
Finally, the fact that the real orientation $x\in
\mur^{\rho}_{\zt}(\cp^{\infty})$ restricts to $1$ under the
isomorphism $\tilde\mur^{2}(S^{\rho})\approx \mur^{0}(\text{pt})$
implies that the element $b_{1}(1)\in H_{1}X^{2}$ corresponds to the
element $[0]\in H_{0}X^{0}$ under the map 
\[
\Sigma X^{0}\to BX_{0}\to X^{2}.
\]
This implies that right column in~\eqref{eq:35} coincides with the
homology suspension induced by
\[
\Sigma X^{2(k-1)}\approx \Sigma \Omega BX^{2(k-1)}\to BX^{2(k-1)}
\]
and so is the edge homomorphism of the Rothenberg-Steenrod spectral
sequence 
\[
\tor_{s,t}^{H_{\ast}X^{2(k-1)}}(\ftwo,\ftwo) \Rightarrow H_{s+t}BX^{2(k-1)}.
\]

\section{Proof of the main theorem}
\label{sec:proof-main-theorem}

Our main result is Theorem~\ref{thm:26} which we restate for the
reader's convenience

\begin{thm}[Theorem~\ref{thm:26}]
\label{thm:34}
The induced map
\begin{equation}
\label{eq:36}
\mot(\hringmurw^{\bullet}) \to \hringmur^{\bullet}
\end{equation}
is an isomorphism.
\end{thm}

As explained in \S\ref{sec:geom-fixed-points} the results of
Ravenel-Wilson~\cite{MR0356030} reduce the above to showing that the
induced map on modified geometric fixed points 
\begin{equation}
\label{eq:37}
\bphig\mot(\hringmurw^{\bullet}) \to \bphig(\hringmur^{\bullet})
\end{equation}
is an isomorphism.
This is a map in $\rmod{H\ftwo}$ which is equivalent to the category
$\grvect$ of graded vector spaces over $\ftwo$.  As explained in
\S\ref{sec:spaces-x2k}, under this equivalence the map
becomes 
\begin{equation}
\label{eq:38}
\Phi^{1}(\hringmurw^{\bullet}) \to H_{\ast}(X^{\bullet}).
\end{equation}
classifying the weight one additive curve $u$ defined in
\S\ref{sec:spaces-x2k}.  We will prove Theorem~\ref{thm:34} by showing
that~\eqref{eq:38} is an isomorphism.  

The universal Hopf ring $\Phi^{1}(\hringmurw^{\bullet})$ will come up
frequently, in many diagrams, and with an additional index indicating
the grading in $\grvect$.  To keep things simple we will use the
abbreviation
\[
\hringr^{\bullet} = \Phi^{1}(\hringmurw^{\bullet}) 
\]
and write~\eqref{eq:38} as
\begin{equation}
\label{eq:39}
\hringr^{\bullet}\to H_{\ast}(X^{\bullet}).
\end{equation}

With this notation we isolate the statement that remains to be proved.

\begin{prop}
\label{thm:35}
The map
\begin{equation}
\label{eq:40}
\hringr^{\bullet}\to H_{\ast}(X^{\bullet}).
\end{equation}
described in \S\ref{sec:spaces-x2k} is an isomorphism.
\end{prop}

Our proof of Proposition~\ref{thm:35} follows the inductive argument
of Chan~\cite{MR675011}.  As mentioned in the introduction, the result
can also be quickly deduced from the more general
result~\cite[Theorem~1.5]{MR2317067} of Kitchloo-Wilson.

\subsection{The stable range}

In spectra there is a canonical weak equivalence 
\[
\varinjlim S^{-2k}\wedge \umu^{2k} \to MU
\]
given by the unit of the adjunction between $\Sigma^{\infty}$ and
$\Omega^{\infty}$, and the map
\[
\varinjlim H_{\bullet}(S^{-2k}\wedge \umu^{2k};\Z) \to H_{\bullet} (MU;\Z)
\]
is an isomorphism.

\begin{lem}
\label{thm:37} 
Let $\cat C$ be an additive category with countable
coproducts.  If 
\[
\mathbf{F}:\Ab\to \cat C
\]
is a coproduct preserving functor, then the colimit
\[
\varinjlim \mathbf{F}(H_{\ast}S^{-2k}\umu^{2k})
\]
exists in $\cat C$ and the map 
\[
\varinjlim \mathbf{F}(H_{\ast}S^{-2k}\umu^{2k})\to \mathbf{F}(H_{\ast}MU)
\]
is an isomorphism.
\end{lem}

\begin{pf}
Since $\mathcal{F}$ preserves coproducts, it suffices to show that the
coequalizer diagram
\begin{equation}
\label{eq:42}
\bigoplus H_{\ast}S^{-2k}\umu^{2k} \rightrightarrows
\bigoplus H_{\ast}S^{-2k}\umu^{2k} \to H_{\ast}MU
\end{equation}
becomes a coequalizer diagram after applying $\mathcal{F}$.   However
since $H_{\ast}(MU;\Z)$ is a free abelian group, the
diagram~\eqref{eq:42} extends to a split coequalizer.   The lemma
follows since split coequalizers are universal colimits.
\end{pf}

Now consider the diagram 
\begin{equation}
\label{eq:43}
\xymatrix@C=2em{
\cdots  \ar[r]  & \mot(H_{\bullet}S^{-2k}\wedge \umu^{2k})
\ar[r]\ar[d]  &\mot(H_{\bullet}S^{-2(k+1)}\wedge \umu^{2(k+1)})
\ar[r]\ar[d]  & \cdots
\\
\cdots   \ar[r]  & H\zm\wedge \umur^{k\rho} \ar[r] & H\zm\wedge
\umur^{(k+1)\rho} \ar[r]& \cdots  \mathrlap{\ .}
}
\end{equation}
By Proposition~\ref{thm:37} with $\mathbf{F}=\mot$, the colimit of the
top row exists in $\hzmod$ and is given by $\mot(H_{\bullet}MU)$.  The
map from the bottom row to $H\zm\wedge \mur$ is a cone on the bottom
row (in the sense of MacLane~\cite[\S III.3, p.~67]{MR1712872}), so there is a unique map
\begin{equation}
\label{eq:44}
\mot(H_{\bullet}MU)\to H\zm\wedge \mur
\end{equation}
of cones on the rows
of~\eqref{eq:43}.   Since the 
map from the top row to the bottom row of~\eqref{eq:43} becomes an isomorphism after
applying $\res$, the map~\eqref{eq:44} is an isomorphism after
applying $\res$.    Since $H\Zm\wedge \mur$ is
pure (Proposition~\ref{thm:22}) this implies that~\eqref{eq:44} is
an isomorphism (Proposition~\ref{thm:18}).  

\begin{prop}
\label{thm:38}
The map $\hringr^{2k}\to H_{\ast}X^{2k}$
is an isomorphism in degrees less than $\max\{1,2k \}$.
\end{prop}

\begin{pf}
In degree $0$, the map $\hringr^{2k}\to H_{\ast}X^{2k}$ is the
identity map of $\ftwo[MU^{2k}]$, so it
suffices to consider the case $k>0$ (in which case $MU^{2k}=0$), and
work in reduced homology.   We therefore wish to show that the
map 
\[
\mot(\tilde H_{\bullet}\umu^{2k}) \to H\zm\wedge \umur^{k\rho}
\]
becomes a $2k$-equivalence\footnote{Specifically, an
epimorphism in degree $2k$ and an isomorphism in lower degrees.}
after applying modified geometric fixed points.

Consider the diagram 
\begin{equation}
\label{eq:45}
\xymatrix{ \mot(\tilde H_{\evenstar}\umu^{2k}) \ar[r]\ar[d] &
\mot(H_{\evenstar}S^{2k}\wedge MU)
\ar[d]^{\sim} \\
 H\zm\wedge \umur^{k\rho}  \ar[r]        & H\zm\wedge
 S^{k\rho}\wedge \mur\mathrlap{\ ,}
}
\end{equation}
extracted from the map of cones on~\eqref{eq:43} by tensoring the
$k^{\text{th}}$-term with $H\zm\wedge S^{k\rho}$.  We wish to show
that the left vertical arrow becomes a $2k$-equivalence after applying
$\bphig$.  By the discussion above, the right vertical arrow is an
isomorphism.  The map
\[
\tilde H_{\evenstar}\umu^{2k}\to H_{\evenstar}S^{2k}\wedge MU
\]
is a $4k$-connected map of evenly graded free abelian groups.   
Since $\bphig\circ\mot = \Phi^{1}$ this
means that the map 
\[
\bphig \mot(\tilde H_{\evenstar}\umu^{2k}) \to
\bphig\mot(H_{\evenstar}S^{2k}\wedge MU)
\]
is a $2k$-connected map of graded $\ftwo$ vector spaces.
For the bottom map, the fact that the spectrum $S^{k\rho}\wedge
\mur$ is equivariantly $(k-1)$-connected implies that the map
\begin{equation}
\label{eq:61}
\Sigma^{\infty}\umur^{k\rho} \to S^{k\rho}\wedge \mur
\end{equation}
is equivariantly $2k$-connected, hence so is the map on geometric fixed
points, since geometric fixed points preserves connectivity.  
The identity
\[
\bphig H\zm\wedge W = H\Z/2\wedge \phig W
\]
then implies that the bottom row of~\eqref{eq:45} becomes a $2k$-equivalence
after passing to modified geometric fixed points.   
It follows that after applying modified geometric fixed points
to~\eqref{eq:45}, the top, right, and bottom arrows become isomorphism
in dimensions less than $2k$.   The same is therefore true of the left
vertical arrow.
\end{pf}

\subsection{Facts about the algebras $\hringr^{2k}$}
\label{sec:facts-about-algebras}

We now collect the facts about the
map~\eqref{eq:39} that will be used in the proof of
Theorem~\ref{thm:34}.

For a graded free module $M$ over a ring $R$ write
\[
\dim_{\alpha}M = \dim_{\alpha}^{R}M=\sum (\dim_{R} M_{k})\alpha^{k}.
\]
Let $\partitions{k}$ be the number of partitions of $k$, defined by the expansion
\[
\prod_{m\ge 1}(1-\alpha^{m})^{-1}= \sum_{k\ge 0} \partitions{k}\alpha^{k},
\]
so that
\begin{align*}
\dim_{\alpha}\pi_{\ast}S^{2k}\wedge MU &= \sum \partitions{m}\alpha^{2(k+m)} \\ 
\dim_{\alpha}\pi_{\ast}(\umu^{2k})' &= \sum_{k+m>0} \partitions{m}\alpha^{2(k+m)}.
\end{align*}
Since  for any connected infinite loop space $Z$ the map 
\[
\sym \pi_{\ast}X\otimes \Q
\to H_{\ast}(X;\Q)
\]
is an isomorphism, and since $H_{\ast}(\umu^{2k})'$ is free abelian, we have
\[
\dim_{\alpha}H_{\ast}(\umu^{2k})' = \prod_{k+m>0}(1-\alpha^{2(k+m)})^{-\partitions{m}}.
\]
It follows from~\ref{thm:5}\footnote{With a little
more investment in setting things up, these purely algebraic facts can
easily be given purely algebraic proofs.  See, for
example~\cite{hill18:_univer_hopf}.} that
\[
\dim_{\alpha}(\hringr^{2k})' = \prod_{k+m>0}(1-\alpha^{k+m})^{-\partitions{m}},
\]
and so 
\begin{equation}
\label{eq:46}
\dim_{\alpha}\rwHringIdeal^{2k}=\frac{\dim_{\alpha}(\hringr^{2k})'}{\dim_{\alpha}H_{\ast}(\umu^{2k})'} 
= \prod_{k+m>0}(1+\alpha^{k+m})^{\partitions{m}}.
\end{equation}
We also know that
\[
\pi_{0}X^{\kk}=\pi_{0}\umu^{2k}=MU^{2k}
\]
is a free abelian group of rank $\partitions{-k}$.

\begin{prop}
\label{thm:39}
The rings $\hringr^{2k}$, $(\hringr^{2k})'$ and the  maps $(\hringr^{2k})'\to H_{\ast}(X^{\kk})'$
have the following properties
\begin{thmList}
\item\label{item:1} The Map $(\hringr^{2k})'\to H_{\ast}(X^{\kk})'$ is an isomorphism in degrees
less than $\max\{1,2k \}$.
\item\label{item:2} The map $b_1(1)\circ(\slot):Q\hringr^{2k}\to \rwHringIdeal^{2(k+1)}\subset 
(\hringr^{2(k+1)})'$ is a monomorphism  and the following diagram commutes
\[
\xymatrix{
Q\hringr^{2k}  \ar[r]\ar[d]_{b_1(1)\circ(\slot)}  & QH_{\ast}X^{\kk} 
\ar[d]^{b_1(1)\circ(\slot)} \\
\hringr^{2(k+1)}  \ar[r]        & H_{\ast}X^{2(k+1)}\mathrlap{\ .}
}
\]
\item\label{item:3}  The algebra $\tor_{s,t}^{\hringr^{2k}}(\ftwo,\ftwo)$ is generated by
elements in $\tor_{1,t}$.   
\item\label{item:4} With respect to the grading on
$\tor_{\ast,\ast}^{\hringr^{2k}}(\ft,\ft))$ in which the term
$\tor_{s,t}$ has degree $s+t$, there is an equality of Poincar\'e
series
\[
\dim_{\alpha}(\tor_{\ast,\ast}^{\hringr^{2k}}(\ft,\ft))
=\dim_{\alpha}(\brk{2(k+1)}).
\]
\end{thmList}
\end{prop}

\begin{pf}
Part~\thmItemref{item:1} is Proposition~\ref{thm:38}.  The first
assertion in part~\thmItemref{item:2} is gotten by applying $\bphig$,
and the commutative diagram is~\eqref{eq:34}.
Part~\thmItemref{item:3} is immediate from the fact that $MU^{2k}$ is
a free abelian group and $(\hringr^{2k})'$ is a polynomial algebra, so
that the $\tor$ algebras in question are exterior algebras on the
elements in degree $\tor_{1,t}$.  For part~\thmItemref{item:4}, we use
that fact that $(\hringr^{2k})'$ is a polynomial algebra and that
\[
\dim_{\alpha}Q (\hringr^{2k})' = \sum_{k+m>0} \partitions{m}\alpha^{k+m},
\]
where $Q$ denotes the indecomposables.   This implies that
\[
\tor_{\ast,\ast}^{(\hringr^{2k})'}(\ft,\ft)
\]
is an exterior algebra on the suspension of the
indecomposables, so that, with respect to total degree,
\[
\dim_{\alpha}\tor_{\ast,\ast}^{(\hringr^{2k})'}(\ft,\ft) = 
\prod_{k+m>0}(1+\alpha^{k+1+m})^{\partitions{m}}.
\]
Since $MU^{2k}$ is a free abelian group of rank $\partitions{-k}$ we
can conclude
\begin{align*}
\dim_{\alpha}\tor_{\ast,\ast}^{\hringr^{2k}}(\ft,\ft) &= 
(1+\alpha)^{\partitions{-k}}\prod_{k+m>0}(1+\alpha^{k+1+m})^{\partitions{m}} \\
&= 
\prod_{k+m+1>0}(1+\alpha^{k+1+m})^{\partitions{m}}.
\end{align*}
The claim now follows from~\eqref{eq:46}.
\end{pf}

\subsection{The induction}
\label{sec:induction}

The assertion of Proposition~\ref{thm:35} is that for all $k$, the
map
\begin{equation}
\label{eq:47}
\hringr^{2k} \to H_{\ast}X^{2k} 
\end{equation}
is an isomorphism.   The argument will involve an induction on degree,
so we are actually interested in the assertion that for some $\ell$,
the map 
\begin{equation}
\label{eq:48}
(\hringr^{2k})_{\ast\le \ell} \to H_{\ast\le \ell}X^{2k} 
\end{equation}
is an isomorphism.  By the discussion in \S\ref{sec:decomposition},
the map~\eqref{eq:47} is the tensor
product of the identity map of $\ft[\pi_{0}X^{2k}]$ with
\begin{equation}
\label{eq:49}
(\hringr^{2k})' \to H_{\ast}(X^{2k})',
\end{equation}
which in turn, is the tensor product of the identity map of
$H_{\ast}(\umu^{2k})'$ with
\begin{equation}
\label{eq:50}
\rwHringIdeal^{2k} \to H_{\ast}BX^{2(k-1)}.
\end{equation}
It follows that the condition that any one of these three maps be an
isomorphism in degrees less than or equal to $\ell$ is equivalent to
the condition that the other two be so.  

Now each of these three maps is a map of Hopf algebras and is an
isomorphism in degrees less than or equal to $\ell$ if and only if the
corresponding map
\begin{equation}
\label{eq:51}
\begin{aligned}
Q\hringr^{2k}&\to QH_{\ast}X^{2k} \\
Q(\hringr^{2k})' &\to QH_{\ast}(X^{2k})' \\
Q\rwHringIdeal^{2k} & \to QH_{\ast}BX^{2(k-1)}
\end{aligned}
\end{equation}
of indecomposables is.   This is clear for the second and third maps
as both sides are connected graded Hopf algebras.   It is so for the
first since the map 
\[
Q\hringr^{2k}\to QH_{\ast}X^{2k}
\]
is the sum of the map
\[
Q(\hringr^{2k})'\to QH_{\ast}(X^{2k})'
\]
and the identity map of 
\[
Q\ft[\pi_{0}X^{2k}] = \pi_{0}X^{2k}\otimes \ft.
\]

Of this, what we will use is
\begin{prop}
\label{thm:40}
For a given $k$ and $\ell$, if any of the following three maps is an
isomorphism
\begin{align*}
(\hringr^{2k})_{\ast\le \ell} &\to H_{\ast\le \ell}(X^{2k}) \\
(\rwHringIdeal^{2k})_{\ast\le \ell} &\to H_{\ast\le \ell}(BX^{2(k-1)}) \\
(Q\hringr^{2k})_{\ast\le \ell} &\to QH_{\ast\le \ell}(X^{2k}) \\
\end{align*}
so are the other two.   \qed
\end{prop}

\begin{figure}
\includegraphics[]{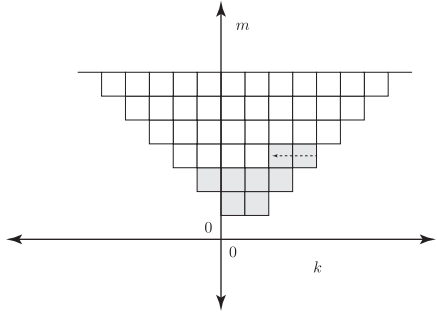}
\caption{The Induction Region}
\label{fig:1}
\end{figure}

Following the inductive the argument of Chan~\cite{MR675011} (and
Ravenel-Wilson~\cite{RW:HR}) consider the assertion that
\begin{equation}
\label{eq:52}
(\hringr^{2k})_{m+k}\to H_{m+k}(X^{2k})
\end{equation}
is an isomorphism.  By part~\thmItemref{item:1} of
Proposition~\ref{thm:39} we know this to be so for $k>m$ ($m+k<2k$)
and for $k< 1-m$ ($m+k<1$).  Among other things this means that it is
an isomorphism for all $m\le 0$.  The region for which we do {\em not}
know it to be an isomorphism is indicated by the boxes in
Figure~\ref{fig:1}.

We will establish, by increasing induction on $m$, that for all $k$
the map~\eqref{eq:52} is an isomorphism.  The induction starts with
$m\le 0$.  

For the inductive step fix $m$ and assume that
\begin{equation}
\label{eq:53}
(\hringr^{2k})_{m'+k}\to H_{m'+k}(X^{2k}).
\end{equation}
is an isomorphism for all $k$ and all $m'<m$.  We will show by
decreasing induction on $k$ that it is an isomorphism for $m'=m$, that
is, that the map~\eqref{eq:52} is an isomorphism.

The induction on $k$ starts with $k=m+1$.  Assume we have
established that~\eqref{eq:52} is an isomorphism for some $k$.   We
wish to conclude that
\begin{equation}
\label{eq:54}
(\hringr^{2(k-1)})_{m+k-1}\to H_{m+k-1}(X^{2(k-1)})
\end{equation}
is an isomorphism.  For this we will use the spectral sequence
\[
\tor_{s,t}^{H_{\ast}(X^{2(k-1)})}(\ft,\ft)
 \Rightarrow H_{s+t}(BX^{2(k-1)}).
\]
to show that the map
\[
(Q\hringr^{2(k-1)})_{m+k-1}\to QH_{m+k-1}(X^{2(k-1)})
\]
is a monomomorphism, and then, by dimension count, an isomorphism.
Combined with the induction hypothesis, this shows that
\[
(Q\hringr^{2(k-1)})_{\ast\le m+k-1}\to QH_{\ast\le m+k-1}(X^{2(k-1)})
\]
is an isomorphism, and so by Proposition~\ref{thm:40} that 
\[
(\hringr^{2(k-1)})_{\ast\le m+k-1}\to H_{\ast\le m+k-1}(X^{2(k-1)})
\]
is as well.   

\begin{figure}
\includegraphics{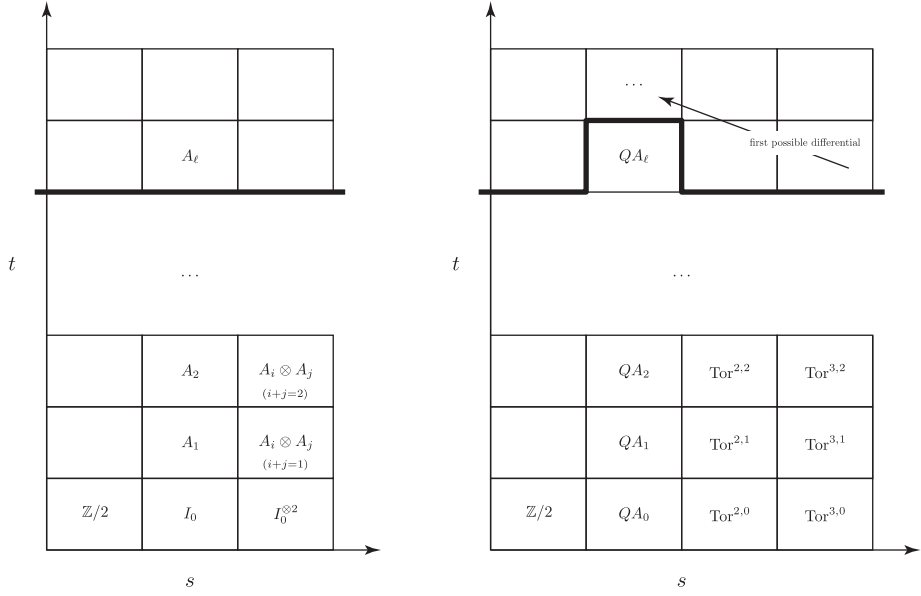}
\caption{The $E_{1}$ and $E_{2}$-terms of the $\tor$ spectral sequence}
\label{fig:2}
\end{figure}

To simplify the notation a bit, write
$A_{\ast}=H_{\ast}(X^{2(k-1)})$, and $I_{\ast}\subset A_{\ast}$ for the
kernel of the counit (augmentation) $A_{\ast}\to \ft$.  Of course, for $t>0$,
$I_{t}=A_{t}$.  The reduced bar construction for computing
\[
\tor_{s,t}=\tor_{s,t}^{A_{\ast}}(\ft,\ft)
\]
is depicted on the left in Figure~\ref{fig:2}.
It serves as an $E^{1}_{s,t}$-term for the spectral sequence, with
$d_{1}$ the bar differential.  The differential $d_{r}$ has bidegree
$(-r,r-1)$.  The $E^{2}_{1,t}$-term consists of the indecomposables
$Q(A)_{t}=(I/I^{2})_{t}$.  All of the elements in the groups
$\tor_{1,t}$ are permanent cycles.  If it happens that for $t<\ell$,
the algebra $\tor_{s,t}^{A_{\ast}}(\ft,\ft)$ is generated by
$\tor_{1,t}$, then there can be no differentials into the terms
$\tor_{s,t}$ for $t\le \ell + 1$.

We apply these considerations with $\ell=m+k-1$.  
Since the map $\hringr^{2(k-1)}\to H_{\ast}(X^{2(k-1)})$ is an isomorphism in
degrees less than $\ell = m+k-1$, we have for $t < \ell$, an isomorphism
\[
\tor_{s,t}^{H_{\ast}(X^{2(k-1)})}(\ft,\ft) \approx
\tor_{s,t}^{\hringr^{2(k-1)}}(\ft,\ft).
\]
Since $\tor_{0,t}=0$ for $t>0$ this holds in particular when
$s+t \le  \ell$ and when $s+t=\ell+1$ provided $s>1$.    By our
induction hypothesis on $m$, 
part~\thmItemref{item:3} and the discussion above, there are no
differentials entering or leaving the $E^{r}_{s,t}$ term when
$s+t=\ell+1$.   This means that the map
\[
b_{1}(1)\circ(\slot):QH_{\ell}(X^{2(k-1)})\to  H_{\ell+1}(BX^{2(k-1)})
\]
is a monomomorphism.  Consider the diagram
\[
\xymatrix@C=3em{
0\ar[r] & (Q\hringr^{2(k-1)})_{\ell}  \ar[rr]^{b_1(1)\circ(\slot)}\ar[d]  &&  (\rwHringIdeal^{2k})_{\ell+1} \ar[d]^{\approx}   &  &\\
0\ar[r] & QH_{\ell}X^{2(k-1)}   \ar[rr]_{b_1(1)\circ(\slot)}        && H_{\ell+1}BX^{2(k-1)}  \ar[r]
& M\ar[r] &0
}
\]
in which the bottom row is exact.  The top map is a monomorphism by
part~\thmItemref{item:2} of Theorem~\ref{thm:39} and the right
vertical map is an isomorphism by our induction hypothesis.  This
establishes the assertion that $(Q\hringr^{2(k-1)})_{\ell}\to
QH_{\ell}(X^{2(k-1)})$ is a monomorphism.

Now for the dimension count.  From the spectral sequence, the vector
space $M$ has a filtration whose associated graded vector space is
\[
\bigoplus_{\substack{s+t = \ell+1 \\ s > 1}}
\tor_{s,t}^{\hringr^{2(k-1)}}(\ft,\ft),
\]
and so 
\begin{equation}
\label{eq:55}
\dim M = \sum_{\substack{s+t = \ell+1 \\ s > 1}}
\dim \tor_{s,t}^{\hringr^{2(k-1)}}(\ft,\ft).
\end{equation}
Now part~\thmItemref{item:4} of Proposition~\ref{thm:39} implies that
\begin{equation}
\label{eq:56}
\begin{aligned}
\dim H_{\ell+1}(BX^{2k}) &=\dim (\rwHringIdeal^{2(k+1)})_{\ell+1} \\
&= \sum_{s+t = \ell+1}\dim
\tor_{s,t}^{\hringr^{2(k-1)}}(\ft,\ft).
\end{aligned}
\end{equation}
The difference between~\eqref{eq:56} and~\eqref{eq:55} is 
\[
\dim \tor^{\hringr^{2(k-1)}}_{1,\ell}(\ft,\ft).
\]
It follows that
\[
\dim QH_{\ell}X^{2(k-1)} =\dim \tor^{\hringr^{2(k-1)}}_{1,\ell}(\ft,\ft)= \dim (Q\hringr^{2(k-1)})_{\ell},
\]
and so the map $(Q\hringr^{2(k-1)})_{\ell}\to QH_{\ell}(X^{2(k-1)})$ must be an
isomorphism.     This completes the proof.

\appendix
\section{Notation and basic notions}
\label{sec:notat-basic-noti}

\subsection*{Strongly additive categories}

An {\em additive category} is a category $\cat C$ satisfying the
following conditions:
\begin{thmList}
\item The category $\cat C$ posses finite products and coproducts.   In
particular $\cat C$ has an initial object and a terminal object,
corresponding to the empty coproduct and product.
\item $\cat C$ is {\em pointed} in the sense that the map from the
initial object to the terminal object is an isomorphism.  This object
is denoted $\ast$ and its existence  makes
all of the hom set $\cat C(X,Y)$ into {\em pointed sets} with base
point the map $X\to \ast \to Y.$

\item For every finite collection of objects $\{X_{i} \}$ of $\cat C$,
the map 
\[
\amalg X_{i} \to \prod X_{i}
\]
whose $i^{\text{th}}$ component is the identity map on the summand $X_{i}$ and the
base point $\ast$ on other summands is an isomorphism.   This equips
the objects $\cat C(X,Y)$ the structure of a commutative monoid.
\item For all $X$ and $Y$, the commutative monoid $\cat C(X,Y)$ is an
abelian group.
\end{thmList}

An additive category closed under arbitry coproducts is called
{\em strongly additive.}

\subsection*{Symmetric monoidal categories}

A {\em  monoidal category} is a
category $\cat C$ equipped with a functor $\otimes:\cat C\times \cat C
\to \cat C$ together with a natural isomorphism
\[
\alpha_{X,Y,Z}:X\otimes (Y\otimes Z) \to (X\otimes Y)\otimes Z 
\]
satisfying the pentagon axiom, and for which there exists an object
$\unit\in \cat C$ which is a two-sided unit for $\otimes$, compatible
with the $\alpha_{X,Y,Z}$ (see~\cite[VII, \S1]{MR1712872}
\cite[Chapter~2]{MR3242743}, and for further
references~\cite{MR1250465}).  In~\cite[Remark~2.2.9]{MR3242743}
Etingof et al.\ point out that the unit, if it exists, is unique up to
unique isomorphism and so corresponds to a property and not a
structure.

A {\em symmetric monoidal category } is a monoidal category $\cat C$
equipped with a natural symmetry isomorphism 
\[
c_{X,Y}:X\otimes Y\to Y\otimes X
\]
for which 
\[
c_{Y,X}\circ c_{X,Y}= \id_{X\otimes Y}
\]
and the ``hexagon'' axiom described in \cite[VII,
\S7]{MR1712872}~\cite[\S8.1]{MR3242743} and~\cite[\S1]{MR1337494}.
Etingof et al.~\cite[Exercise~8.1.6]{MR3242743} point out that the
identities~\cite[VII \S7 (2)]{MR1712872} are a consequence of the
other axioms.

The functors between symmetric monoidal categories we consider are the
{\em strong symmetric monoidal functors} in the sense of
Thomason~\cite[\S1]{MR1337494}.  (See also~\cite[\S2.4]{MR3242743} for
an incisive and full discussion.)

We use $(\cat C, \otimes)$ to indicate the data of a symmetric
monoidal category, rather than $(\cat C,\otimes, \alpha, c)$, and denote the
unit by $\unit$.

In this paper a
{\em strongly
additive} symmetric monoidal category is a symmetric monoidal category
$(\cat C,\otimes)$ in which $\cat C$ is additive, contains
arbitrary coproducts, and the monoidal functor
$\otimes:\cat C\times\cat C\to \cat C$ commutes with coproducts in each variable.  
In particular the monoidal product is biadditive.   A {\em strongly additive
symmetric monoidal functor} between strongly additive symmetric monoidal
categories is a (strong) symmetric monoidal functor which commutes
with coproducts.  

\subsection*{Grading}
In this paper we are frequently
concerned with evenly graded abelian groups, and it is useful to use
the symbol $\evenstar$ as a wildcard matching even integers.  So for
example, for a space or spectrum $X$, the symbols $H_{\evenstar}X$ and
$\pi_{\evenstar}X$ denote the even dimensional homology and
homotopy groups of $X$ regarded as evenly graded abelian groups.
Similarly an evenly graded Hopf ring is denoted by something like
$\hringh^{\evenstar}$ with $\hringr^{n}$ indicating the component of degree $n$ when
$n$ is even, and undefined if $n$ is odd.

In equivariant homotopy theory we follow the convention,
introduced by Hu and Kriz in~\cite{MR1808224}, of using the symbol
$\star$ as a wildcard matching a real virtual representation, and
using $M_{\star}$ and $M^{\star}$ to denote $RO(G)$-graded objects.

Given an $RO(G)$-graded object $M^{\star}$ and a representation
$V$ of $G$ of dimension $d$, there is a $d\Z$-graded object $M(V)$ with
$M(V)^{nd}=M^{nV}$ (\S\ref{sec:grading}).

For ordinary $\Z$-grading we use, as is customary, the symbol
$\ast$ as in $M_{\ast}$ and $M^{\ast}$.

\subsection*{Yoneda}

We do not distinguish in notation between an object $X$ in a
category $\cat D$ and the functor represented by $X$.   Thus for  $X,Y\in
\cat D$ we write $X(Y)=\cat D(Y,X)$.   This comes up most of the time
when $\cat D=\coalg\cat C$ is the category of coalgebras in a
symmetric monoidal category $\cat C$.

\subsection*{Spaces in a spectrum}

If $E$ is a spectrum and $n\in \Z$ we write
\[
\underline{E}^{n}=\Omega^{\infty}\Sigma^{n}E
\]
for the $n^{\text{th}}$ space in the associated $\Omega$-spectrum.
This is a departure from the long established convention found
throughout the literature%
~\cite{Lima,MR0116332,MR0107862,Whitehead:ght,
Ad:SHGH,adams69:_stabl,LMayS,
WSWThesis}.
It is
with some reluctance that we have chosen to make this convention, but
it leads to clearer expressions in the situations of this paper.  For
example the $\underline{E}^{n}$ represents the fuctor $E^{n}(X)$, so
this is consistent with our convention for the Yoneda embedding in the
case $\cat C=\ho\spaces$ is the homotopy category of spaces
\[
E^{n}(X)=\underline{E}^{n}(X)=\ho\spaces(X,\underline{E}^{n}).
\]

\subsection*{Further miscellaneous notation}

\begin{itemize}

\item $\coalg\cat C$ is the category of cocommutative counital
coalgebras in a symmetric monoidal category $\cat C$ (\S\ref{sec:hopf-rings-1}).

\item $\hopfrings_{\rbullet{R}}(\cat C)$ is the category of evenly
graded Hopf ring in $\cat C$ over the evenly graded ring $\rbullet{R}$.

\item  From \S\ref{sec:universal-hopf-rings} to
\S\ref{sec:hopf-ring-xbullet} the symbol $H_{\ast}X$ refers the
homology $H_{\ast}(X;\Z)$ of a space or spectrum $X$.  From
\S\ref{sec:hopf-ring-xbullet} until the end of the paper
$H_{\ast}(X)$ denotes homology with coefficients in $\ftwo$.

\item For a commutative ring $k$, the category $\gmodules{k}$ is the
symmetric monoidal category of left $k$-modules with monoidal product
the tensor product and symmetry constraint given by the Koszul sign
rule (Example~\ref{eg:3}).  The category $\gmodules{k}$ is a strongly
additive symmetric monoidal category.  When $k=\Z$ we write $\abstar$
instead of $\gmodules{\Z}$ (\S\ref{sec:ravenel-wilson-hopf-2}).

\item For an (equivariant) $\einfty$ ring spectrum $R$, the category
$\rmod{R}$ is the homotopy category of (equivariant) left $R$-modules
(Example~\ref{eg:4} and \S\ref{sec:equiv-r-modul}).  The category
$\rmod{R}$ is a strongly additive symmetric monoidal category under
$X\underset{R}{\wedge}Y$.

\item For a ($G$-)space $X$, the free (equivariant) $R$-module
generated by $X$ is $\freeRmodule{R}{X}=R\wedge
X_{+}$ (Example~\ref{eg:4} and \S\ref{sec:equiv-r-modul}).

\item $\A^{1}$ is the Hopf algebra $H_{\ast}(\cp^{\infty})$
(\S\ref{sec:ravenel-wilson-hopf-2}).  It has a graded basis
$\{\beta_{i}\mid i=0,1,\cdots\}$ in which $\beta_{i}$ is the generator
of $H_{2i}(\cp^{\infty};\Z)$ dual to $x^{i}$, where $x\in
H^{2}(\cp^{\infty};\Z)$ is the first Chern class of the tautological
line bundle.

\item $\hringmurw^{\bullet}$ is the Ravenel-Wilson Hopf ring (Definition~\ref{def:6}).

\item $\cat D^{\bullet}$ denotes the category of evenly graded objects
of a category $\cat D$ (\S\ref{sec:evenly-grade-hopf}).

\item $\frab$ is the strongly additive symmetric monoidal category of
evenly graded free abelian groups and tensor product (\S\ref{sec:gener-raven-wils}).

\item $\motZ:\frab\to \zmod$~\eqref{eq:7} is the unique strongly additive functor satisfying 
$\motZ(\Z[2n])= H\Z\wedge S^{2n}$.

\item $\grvect=\gmodules{\ftwo}$ is the (strongly additive) symmetric
monoidal category of graded vectors spaces over $\ftwo$ and tensor
products (\S\ref{sec:weight-k-curves})

\item $\vmot:\grvect\to \hztmod$~\eqref{eq:7} is the unique strongly additive functor satisfying 
$\vmot(\ft[k])= \hft\wedge S^{2k}$ (\S\ref{sec:weight-k-curves}).

\item if $X$ is a pointed space then $X'$ denotes the connected
component of the base point.

\item If $C$ is a $(-1)$-connected coalgebra in $\gmodules{k}$
equipped with a coalgebra map $b:k\to A_{0}$, then $C'$ denotes the
connected component of $C$ containing $b$ (Definition~\ref{def:11}).   

\item For a coalgebra $C$ the set $\pi_{0}C=C(\unit)$ is the set of
coalgebra maps $\unit\to C$ (Definition~\ref{def:17}).

\item $\Phi^{k}(\slot):\frab\to \grvect$ is the strongly additive
symmetric monoidal functor sending $\Z[2n]$ to $\ftwo[k n]$ (\S\ref{sec:weight-k-curves}).

\item $(\slot)^{\phi}:\grvect\to\grvect$ is the functor which doubles
the grading: $V^{\phi}_{2n}=V_{n}$ (\S\ref{sec:verschiebung}).

\item $\versch:C\to C^{\phi}$ is the Verschiebung
(\S\ref{sec:verschiebung}).   

\item $\A^{1}(k)$ is the Hopf algebra $\Phi^{k}\A^{1}\in \grvect$, and has a
graded basis $\{\beta_{i}(k)\mid i=0,1,\dots \}$ with
$\beta_{i}(k)=\Phi^{k}\beta_{i}$ (Example~\ref{eg:18}).

\item $Q\mathcal H$ is the module of indecomposables in a Hopf algebra
$\mathcal H$~\eqref{eq:63}.

\item $\rho$ is the real regular representation of $\zt$.

\item $\mot:\frab\to \hzmodf$ the unique strongly additive functor satisfying 
$\mot(\Z[2n])= S^{n\rho}\wedge H\zm$ (\S\ref{sec:pure-modules-graded}).

\item Given an $RO(G)$-graded object $M^{\star}$ and a representation
$V$ of $G$ of dimension $d$, $M(V)$ is the $d\Z$-graded object with
$M(V)^{nd}=M^{nV}$ (\S\ref{sec:grading}).

\item $\bphig$ is the modified geometric fixed point functor (\S\ref{sec:geom-fixed-points}).

\item $\hringr^{\bullet}=\Phi^{1}\hringmurw^{\bullet}$ (\S\ref{sec:proof-main-theorem})
\end{itemize}


\def\cprime{$'$}
\providecommand{\bysame}{\leavevmode\hbox to3em{\hrulefill}\thinspace}
\providecommand{\MR}{\relax\ifhmode\unskip\space\fi MR }
\providecommand{\MRhref}[2]{%
  \href{http://www.ams.org/mathscinet-getitem?mr=#1}{#2}
}
\providecommand{\href}[2]{#2}

\end{document}